\documentclass[12pt]{amsart}
\usepackage{amssymb,latexsym}
\usepackage{array}
\usepackage{tikz-cd}
\usepackage{amsmath, amssymb, amsfonts, amstext, amsthm, amscd}
\usepackage[mathscr]{euscript}
\usepackage{enumerate}
\usepackage{tikz,color,soul}
\usepackage[Symbolsmallscale]{upgreek}
\usetikzlibrary{arrows.meta, positioning,arrows}
\usetikzlibrary{matrix,arrows}

\usepackage[english]{babel}
\usepackage{chngcntr}
\usepackage{textgreek}
\newcommand{\dequalto}{\mathrel{
\rotatebox[origin=c]{45}{$=$}}
}

\begin{document}

\thispagestyle{empty} 

\title{Equivariant $K$-theory of cellular toric varieties}

\author{V. Uma}
\address{Department of Mathematics, Indian Institute of Technology Madras, Chennai-600036}
\email{vuma@iitm.ac.in}

\subjclass[2020]{19L47, 14M25}


\keywords{Equivariant K-theory, cellular toric varieties}

\date{}
\begin{abstract}
  In this article we describe the $T_{comp}$-equivariant topological
  $K$-ring of a complete $T$-{\it cellular} toric variety. We further
  show that $K_{T_{comp}}^0(X)$ is isomorphic as an
  $R(T_{comp})$-algebra to the ring of piecewise Laurent polynomial
  functions on the associated fan denoted $PLP(\Delta)$. Furthermore,
  we compute a basis for $K_{T_{comp}}^0(X)$ as a $R(T_{comp})$-module
  and multiplicative structure constants with respect to this
  basis.\end{abstract}

\maketitle

\def\theequation
  {\arabic{section}.\arabic{equation}}

\newcommand{\codim}{\mbox{{\rm codim}$\,$}}
\newcommand{\stab}{\mbox{{\rm stab}$\,$}}
\newcommand{\lr}{\mbox{$\rightarrow$}}

\newcommand{\be}{\begin{equation}}
\newcommand{\ee}{\end{equation}}

\newtheorem{guess}{Theorem}[section]
\newcommand{\bth}{\begin{guess}$\!\!\!${\bf }~}
\newcommand{\eeth}{\end{guess}}
\renewcommand{\bar}{\overline}
\newtheorem{propo}[guess]{Proposition}
\newcommand{\bpropo}{\begin{propo}$\!\!\!${\bf }~}
\newcommand{\epropo}{\end{propo}}

\newtheorem{assum}[guess]{Assumption}
\newcommand{\bas}{\begin{assum}$\!\!\!${\bf }~}
\newcommand{\eas}{\end{assum}}

\newtheorem{lema}[guess]{Lemma}
\newcommand{\blem}{\begin{lema}$\!\!\!${\bf }~}
\newcommand{\elem}{\end{lema}}

\newtheorem{defe}[guess]{Definition}
\newcommand{\bdefe}{\begin{defe}$\!\!\!${\bf }~}
\newcommand{\edefe}{\end{defe}}

\newtheorem{coro}[guess]{Corollary}
\newcommand{\bcor}{\begin{coro}$\!\!\!${\bf }~}
\newcommand{\ecor}{\end{coro}}

\newtheorem{rema}[guess]{Remark}
\newcommand{\brem}{\begin{rema}$\!\!\!${\bf }~\rm}
\newcommand{\erem}{\end{rema}}

\newtheorem{exam}[guess]{Example}
\newcommand{\beg}{\begin{exam}$\!\!\!${\bf }~\rm}
\newcommand{\eeg}{\end{exam}}

\newtheorem{notn}[guess]{Notation}
\newcommand{\bnot}{\begin{notn}$\!\!\!${\bf }~\rm}
\newcommand{\enot}{\end{notn}}

\newcommand{\ch}{{\mathcal H}}
\newcommand{\cf}{{\mathcal F}}
\newcommand{\cd}{{\mathcal D}}
\newcommand{\cR}{{\mathcal R}}
\newcommand{\cv}{{\mathcal V}}
\newcommand{\cn}{{\mathcal N}}
\newcommand{\lra}{\rightarrow}
\newcommand{\ra}{\rightarrow}
\newcommand{\blr}{\Big \rightarrow}
\newcommand{\da}{\Big \downarrow}
\newcommand{\ua}{\Big \uparrow}
\newcommand{\hra}{\mbox{{$\hookrightarrow$}}}
\newcommand{\rt}{\mbox{\Large{$\rightarrowtail$}}}
\newcommand{\dua}{\begin{array}[t]{c}
\Big\uparrow \\ [-4mm]
\scriptscriptstyle \wedge \end{array}}
\newcommand{\ctext}[1]{\makebox(0,0){#1}}
\setlength{\unitlength}{0.1mm}
\newcommand{\cl}{{\mathcal L}}
\newcommand{\cp}{{\mathcal P}}
\newcommand{\ci}{{\mathcal I}}
\newcommand{\bz}{\mathbb{Z}}
\newcommand{\cs}{{\mathcal s}}
\newcommand{\ce}{{\mathcal E}}
\newcommand{\ck}{{\mathcal K}}
\newcommand{\cz}{{\mathcal Z}}
\newcommand{\cg}{{\mathcal G}}
\newcommand{\cj}{{\mathcal J}}
\newcommand{\cc}{{\mathcal C}}
\newcommand{\ca}{{\mathcal A}}
\newcommand{\cb}{{\mathcal B}}
\newcommand{\cx}{{\mathcal X}}
\newcommand{\co}{{\mathcal O}}
\newcommand{\bq}{\mathbb{Q}}
\newcommand{\bt}{\mathbb{T}}
\newcommand{\bh}{\mathbb{H}}
\newcommand{\br}{\mathbb{R}}
\newcommand{\bl}{\mathbf{L}}
\newcommand{\wt}{\widetilde}
\newcommand{\im}{{\rm Im}\,}
\newcommand{\bc}{\mathbb{C}}
\newcommand{\bp}{\mathbb{P}}
\newcommand{\ba}{\mathbb{A}}
\newcommand{\spin}{{\rm Spin}\,}
\newcommand{\ds}{\displaystyle}
\newcommand{\tor}{{\rm Tor}\,}
\newcommand{\bff}{{\bf F}}
\newcommand{\bs}{\mathbb{S}}
\def\ns{\mathop{\lr}}
\def\nssup{\mathop{\lr\,sup}}
\def\nsinf{\mathop{\lr\,inf}}
\renewcommand{\phi}{\varphi}
\newcommand{\tT}{{\widetilde{T}}}
\newcommand{\tG}{{\widetilde{G}}}
\newcommand{\tB}{{\widetilde{B}}}
\newcommand{\tC}{{\widetilde{C}}}
\newcommand{\tW}{{\widetilde{W}}}
\newcommand{\tphi}{{\widetilde{\Phi}}}
\noindent

\section{Introduction}\label{Introduction}

We shall consider varieties over the field of complex numbers unless
otherwise specified.

Let $T\simeq (\mathbb{C}^*)^n$ denote the complex algebraic torus of
dimension $n$.  A {\em $T$-cellular variety} is a normal complex
algebraic $T$-variety $X$ equipped with a $T$-stable algebraic cell
decomposition. In other words there is a filtration \be\label{filter}
X=Z_1\supseteq Z_2\supseteq\cdots\supseteq Z_m\supseteq
Z_{m+1}=\emptyset\ee where each $Z_i$ is a closed $T$-stable
subvariety of $X$ and $Z_i\setminus Z_{i+1}=Y_i$ is $T$-equivariantly
isomorphic to the affine space ${\mathbb{C}}^{n_i}$ equipped with a
linear action of $T$ for $1\leq i\leq m$. Furthermore, $Y_i$ for
$1\leq i\leq m$, are the Bialynicki-Birula cells associated to a {\it
  generic} one-parameter subgroup of $T$ (see Definition
\ref{cellular}). A one-parameter subgroup $\lambda:\mathbb{G}_m\lra T$
is said to be {\em generic} if
$X^{\lambda}=X^{T}=\{ x_1,\ldots, x_m\}$. We note here that the
Bialynicki-Birula cell decomposition of $X$ as well as the property of
$X$ being $T$-cellular both depend on the choice of $\lambda$ (see
Section \ref{cellular varieties}).


Any smooth projective complex algebraic variety $X$ with $T$-action
having only finitely many $T$-fixed points is $T$-cellular.  Indeed
for any projective variety the Bialynicki-Birula cell decomposition is
filtrable (see \cite{BB1}) and since $X$ is smooth the
Bialynicki-Birula cells are smooth (see \cite[Section 3.1]{Br2}).

Let $X=X(\Delta)$ be a toric variety associated to a fan $\Delta$ in
the lattice $N=\bz^n$ with the dense torus $T$. We further assume that
$X(\Delta)$ is $T$-{\it cellular}. We shall give a combinatorial
characterization on the fan $\Delta$ so that the associated toric
variety $X(\Delta)$ is $T$-cellular (see Theorem \ref{combcell}) and
we call such a fan a {\em cellular fan} (see Definition \ref{cellular
  fan}).

Our definition of a cellular variety and a cellular toric variety was
motivated by the definition of a divisive weighted projective space
due to Harada, Holm, Ray, and Williams (see \cite{hhrw}) and that of a
retractable toric orbifold in \cite{saruma}. We were also motivated by
the definition of (possibly singular) varieties with {\it good
  decompositions} due to Carrell and Goresky (see \cite{cago} and
\cite[Section 4.3, Theorem 4.13]{car}).  The $T$-cellular varieties
which we consider are varieties which admit good decompositions in the
above sense.


Examples of $T$-cellular toric varieties include smooth projective
toric varieties \cite{f}, \cite{o}, smooth semi-projective toric
varieties defined in \cite[Section 2]{hs} and the divisive weighted
projective spaces defined in \cite{hhrw}. The retractable toric
orbifolds considered in \cite{saruma} includes projective simplicial
toric varieties such that the $T$-equivariant rational Bialynicki
Birula cells in the cellular decomposition of $X(\Delta)$ are in fact
smooth integral cells (see Remark \ref{rational cellular},
\cite[pp. 101-105]{f}, \cite{Br4}, \cite[Section 3]{gon1}). There are
smooth complete non-projective toric varieties which are $T$-cellular
as seen in the example in \cite[p. 71]{f}. We give an example of a
simplicial non-smooth and non-complete fan (hence the corresponding
toric variety) which is cellular (see Example \ref{example1}). We also
give an example of a non-simplicial, complete fan (hence the
corresponding toric variety) which is cellular (see Example
\ref{non-simplicial cellular fan}).

Let $T_{comp}\simeq (S^1)^n$ denote the maximal compact subgroup of
$T$. We consider the natural restricted action of $T_{comp}\subset T$
on $X(\Delta)$. Let $K^0$ denote the topological $K$-ring (see
\cite{at}) and for a compact connected Lie group $G_{comp}$,
$K^0_{G_{comp}}$ denotes $G_{comp}$-equivariant $K$-ring. In particular,
$K^0_{T_{comp}}$ denotes $T_{comp}$-equivariant topological $K$-ring
(see \cite{segal}).

\subsection{Preliminaries and earlier related results}
Let $X$ be a compact $G_{comp}$-space for a compact Lie group
$G_{comp}$. By $K^0_{G_{comp}}(X)$ we mean the Grothendieck ring of
$G_{comp}$-equivariant topological vector bundles on $X$ with the
abelian group structure given by the direct sum and the multiplication
given by the tensor product of equivariant vector bundles. In
particular, $K^0_{G_{comp}}(pt)$, where $G_{comp}$ acts trivially on
$pt$, is the Grothendieck ring $R(G_{comp})$ of complex
representations of $G_{comp}$. The ring $K^0_{G_{comp}}(X)$ has the
structure of $R(G_{comp})$-algebra via the map
$R(G_{comp})\ra K_{G_{comp}}^0(X)$ which takes
$[V]\mapsto \mathbf{V}$, where $\mathbf{V}=X\times V$ is the trivial
$G_{comp}$-equivariant vector bundle on $X$ corresponding to the
$G_{comp}$-representation $V$. Let $pt$ be a $G_{comp}$-fixed point of
$X$ then the reduced equivariant $K$-ring
${\widetilde K}_{G_{comp}}^0(X)$ is the kernel of the map
$K_{G_{comp}}^0(X)\ra K_{G_{comp}}^0(pt)$, induced by the restriction
of $G_{comp}$-equivariant vector bundles. For $n\in \mathbb{N}$, we
define
$\widetilde{K}^{-n}_{G_{comp}}(X):=\widetilde{K}^0_{G_{comp}}(S^{n}X)$
where $S^{n}X$ is the $n$-fold reduced suspension of $X$.

If $X$ is a locally compact Hausdorff space but not compact we shall
consider $G_{comp}$-equivariant $K$-theory with compact support
denoted by $K_{G_{comp},c}^0(X)$. This can be identified with
${\widetilde K}_{G_{comp}}^0(X^+)$ where $X^+$ denotes the one point
compactification with the assumption that the point at infinity is
$G_{comp}$-fixed. We define
${K}_{G_{comp},c}^{-n}(X):=\widetilde{K}^{-n}_{G_{comp}}(X^+)$. When
$X$ is already compact, we define $X^+=X\sqcup {pt}$ as the disjoint
union of $X$ and a base point. In this case we see that
$K^{0}_{G_{comp},c}(X)=\widetilde{K}^0_{G_{comp}}(X^+)=K^{0}_{G_{comp}}(X)$
since
$K^0_{G_{comp}}(X^+)=K_{G_{comp}}^0(X)\oplus K_{G_{comp}}^0(pt)$. In
\cite[Definition 2.8]{segal}, $K^{0}_{G_{comp},c}$ is denoted by
$K^{0}_{G_{comp}}$ without the subscript $c$. We note here that the
equivariant $K$-theory with compact support is not a homotopy
invariant unlike $K_{G_{comp}}^0$ (see \cite[Proposition 2.3]{segal}).
For example
$K^{0}_{G_{comp}, c}(\mathbb{R}^1)={\widetilde
  K}_{G_{comp}}^0(S^1)=K_{G_{comp}}^{-1}(pt)=0$ whereas
$K^0_{G_{comp},c}(pt)=K^0_{G_{comp}}(pt)=R(G_{comp})$.


We have the equivariant Bott periodicity
$K_{G_{comp}}^{-n}(X)\simeq K_{G_{comp}}^{-n-2}(X)$ given via
multiplication by the Bott element in $K_{G_{comp}}^{-2}(pt)$. (See
\cite{segal}, \cite{at} and \cite{fo}.) This enables us to define
$K^{n}_{G_{comp}}(X)$ for a positive $n\in \mathbb{Z}$ as
$K^{n-2q}(X)$ for $q\geq n/2$.

More generally, we can consider the $G_{comp}$-equivariant generalized
complex oriented cohomology theory (see \cite{ttd} or \cite{mat1}) for
compact $G_{comp}$-CW complexes. We denote this by
$\mbox{h}^*_{G_{comp}}$. If $X$ is a non-compact $G_{comp}$-CW complex
then we can similarly define generalized equivariant cohomology with
compact support by
$\mbox{h}_{G_{comp},c}^0(X)={\widetilde{\mbox{h}}}^0_{G_{comp}}(X^+)$.

Let $X$ be an algebraic variety with the action of an algebraic group
$G$. Then $\mathcal{K}^0_{G}(X)$ denotes the Grothendieck ring of
equivariant algebraic vector bundles on $X$ and $\mathcal{K}_0^{G}(X)$
the Grothendieck group of equivariant coherent sheaves on $X$ (see
\cite{Th}, \cite{mer}). The natural map
$\mathcal{K}^0_{G}(X)\ra \mathcal{K}_0^G(X)$ obtained by sending a
class of a $G$-equivariant vector bundle $\mathcal{V}$ on $X$ to the
dual of its sheaf of local sections is an isomorphism when $X$ is
smooth, but not in general. Moreover, when $G$ is a complex reductive
algebraic group and $G_{comp}$ is a maximal compact subgroup of $G$,
then any complex algebraic $G$-variety $X$ is a $G_{comp}$-space. When
$X$ is a smooth compact complex algebraic $G$-variety, we have a
natural map $\mathcal{K}_0^{G}(X)\ra K^0_{G_{comp}}(X)$ obtained by
first identifying $\mathcal{K}^{G}_0(X)$ with $\mathcal{K}_{G}^0(X)$
and then viewing an algebraic $G$-vector bundle as a topological
$G_{comp}$-vector bundle on $X$ (see \cite[Section 5.5.5]{cg}).  We
note here that a comparison of the topological equivariant $K$-theory
with the algebraic equivariant $K$-theory was analysed in
\cite[Section A.3]{frw}.

Recall that in \cite[Section 5.5]{cg} we have an alternate notion of
$G$-cellular variety $X$, where $X$ admits a decreasing filtration
\eqref{filter} by $G$-stable closed subvarieties $Z_i$ such that
$Y_i=Z_i\setminus Z_{i+1}$ are complex affine spaces equipped with a
linear $G$-action. However, $Y_i$ are not assumed to be the cells of a
Bialynicki-Birula decomposition of $X$. Thus for $T$-varieties, this
notion of cellular is weaker than our definition (see Definition
\ref{cellular}). Recall that in \cite{spal}, Spaltenstein has shown the
existence of affine cell decomposition for Springer fibers of type
$A$. However in \cite[Section 3, Example (3)]{cg}, Carrell and Goresky
show that the affine decomposition of \cite{spal} arises as a plus
cell decomposition associated to a one-parameter subgroup with
isolated fixed points, only in special cases. It follows from
\cite[p. 272]{cg}, that the map
$\mathcal{K}_0^{G}(X)\ra K^0_{G_{comp}}(X)$ is an isomorphism when $X$
is smooth and $G$-cellular. Also see the results in \cite{Bros} and
\cite{k}.

For singular $G$-varieties there are no natural isomorphisms
$\mathcal{K}^0_{G}(X)\cong \mathcal{K}_0^G(X)$ and
$\mathcal{K}_0^{G}(X)\cong K^0_{G_{comp}}(X)$.  On the other hand, for
a singular variety $X$, there is an analogous notion of algebraic
operational equivariant $K$-theory denoted
$\mbox{op} ~\mathcal{K}^0_{G}(X)$, which has been studied and
described for arbitrary toric varieties as piecewise Laurent
polynomial functions on the associated fan (see \cite[Theorem
1.6]{ap}). The results of Anderson and Payne \cite{ap} on algebraic
operational $T$-equivariant $K$-ring of toric varieties generalize the
results of Vezzosi and Vistoli \cite{vv} on the algebraic equivariant
$K$-ring of smooth $T$-toric varieties. In \cite{vv}, Vezzosi and
Vistoli prove that the algebraic $T$-equivariant $K$-ring of a smooth
toric variety is isomorphic to the Stanley-Reisner ring of the
associated fan. (Also see \cite{su} and \cite{s} on the presentation
of ordinary algebraic and topological $K$-theory of smooth projective
toric varieties and toric bundles and smooth complete toric varieties
using alternate techniques.) For singular varieties, the Grothendieck
ring of algebraic vector bundles need not be finitely generated and
can be quite large compared to the Grothendieck group of coherent
sheaves (see \cite{gu}). We remark here that when $G=1$ and $X$ the
elliptic curve, then the Grothendieck group of coherent sheaves as
well as the Grothendieck group of algebraic vector bundles are both
infinitely generated and are also isomorphic. However, the algebraic
operational equivariant $K$-ring seems to behave better in relation to
the Grothendieck group of equivariant coherent sheaves (see
\cite[Theorem 1.3]{ap}) and also the topological equivariant $K$-ring
(we are unable to find a precise reference for the latter assertion
which is probably known).

In the topological direction, in \cite{hhrw}, Harada, Holm, Ray and
Williams study the topological equivariant $K$-ring of weighted
projective spaces $\mathbb{P}(\chi_0,\chi_1,\ldots, \chi_n)$ such that
$\chi_{j-1}$ divides $\chi_{j}$ for $1\leq j\leq n$. They called such
weighted projective spaces {\it divisive} and showed that they have a
$T_{comp}$-invariant cellular structure (see \cite[Section
2]{hhrw}). In \cite[Theorem 5.5]{hhrw} they show that the
$T_{comp}$-equivariant topological $K$-ring of a divisive weighted
projective space is isomorphic to the ring of piecewice Laurent
polynomial functions on the associated fan. In the non-equivariant
setting, the results on the $K$-theory of weighted projective spaces
are due to Al Amrani (see \cite{A1, A2}).

In \cite{saruma}, Sarkar and the author study the topological
equivariant $K$-ring of a toric orbifold which has an invariant
cellular structure. Such a toric orbifold was called {\it retractable}
since the sufficient condition for the toric orbifold to have an
invariant cellular structure was given by the notion of {\it a
  retraction sequence} of the associated simple poytope (also see the
results on the integral cohomology of toric orbifolds by Bahri, Sarkar
and Song \cite{bss}). Indeed the simplicial projective toric varieties
which correspond to simplicial polytopal fans are examples of toric
orbifolds. Therefore, by the results in \cite[Section 4, Theorem
4.2]{saruma}, it follows that for a cellular simplicial projective
toric variety the topological $T_{comp}$-equivariant $K$-ring is
isomorphic to the $R(T_{comp})$-algebra of piecewise Laurent
polynomial functions on the associated fan. This description agrees
with the results of Anderson and Payne on the algebraic equivariant
operational $K$-ring. The corresponding result for topological
$T_{comp}$-equivariant $K$-ring of a smooth projective toric variety
was proved in \cite[Theorem 7.1, Corollary 7.2]{hhrw}.

\subsection{Outline of our main results}

Our aim in this paper is to study the topological equivariant $K$-ring
for certain classes of (possibly singular) toric varieties which are
$T$-cellular (see Definition \ref{cellular}).

Let $X=X(\Delta)$ be a complete $T$-cellular toric variety associated
to the fan $\Delta$.  We first give a GKM type description of the
$T_{comp}$-equivariant topological $K$-ring of $X(\Delta)$ (see
Theorem \ref{main1}), using methods similar to those used by Mc Leod
in \cite[Section 1.5, Theorem 1.6]{ml} for the description of the
topological $T_{comp}$-equivariant $K$-ring of a flag variety.

In Section \ref{pl}, we further show that for a complete $T$-cellular
toric variety $X=X(\Delta)$, $K^0_{T_{comp}}(X)$ is isomorphic as an
$R(T_{comp})=R(T)$-algebra to the ring of piecewise Laurent polynomial
functions on $\Delta$ which we denote by $PLP(\Delta)$ (see Theorem
\ref{main2}). In particular, this extends the results in \cite[Theorem
5.5]{hhrw} for a divisive weighted projective space and in
\cite[Theorem 4.2]{saruma} for a cellular simplicial projective toric
variety to any complete cellular toric variety.

By \cite{vv}, for a smooth cellular $T$-toric variety $X(\Delta)$,
$\mathcal{K}^0_{T}(X(\Delta))\cong K^0_{T_{comp}}(X(\Delta))$ is
isomorphic to the Stanley-Reisner ring of $\Delta$. Our results are
therefore a generalization of the results of Vezzosi and Vistoli
\cite{vv}, in the direction of topological equivariant $K$-theory, to
the class of complete cellular (possibly non-simplicial) toric
varieties.

In Section \ref{basis1}, for $X$ a $T$-cellular complete toric
variety, we determine a basis for $K^0_{T_{comp}}(X)$ as a
$R(T_{comp})$-module (see Theorem \ref{basis}) (see \cite[Proposition
4.3]{haheho}, \cite[Proposition 2.13]{ty} and \cite[Section 6]{gon1}
for similar results on the equivariant cohomology ring).  Furthermore,
in Section \ref{strconst}, we determine the multiplicative structure
constants with respect to this basis (see Theorem \ref{coefficients},
Corollary \ref{multstrconst}).

\subsection{Statements of our main results}

Let $X=X(\Delta)$ be a complete $T$-cellular toric variety.  Let
$\Delta(n)=\{\sigma_1,\ldots, \sigma_n\}$ be the set of maximal
dimensional cones in $\Delta$.

Let $x_i$ denote the $T$-fixed point of $X$ corresponding to
$\sigma_i$ for $1\leq i\leq m$. Let $X(\Delta)$ be $T$-cellular with
respect to a generic one-parameter subgroup $\lambda_v$ for $v\in
N$. Let $x_1<\cdots <x_m$ be the ordering of the $T$-fixed points of
$X(\Delta)$ induced by the filtrable Bialynicki-Birula cell
decomposition with respect to $\lambda_{v}$. Let
$Y_i\cong \mathbb{C}^{n-k_i}$ be the plus Bialynicki-Birula cell
correponding to $\lambda_v$. As a $T$-representation $Y_{i}$ is
isomorphic to the direct sum of one-dimensional representations
corresponding to the characters
$u_{i_j}\in M=\mbox{Hom}(N,\mathbb{Z})$ for $1\leq j\leq n-k_i$.


Let $\mathcal{A}$ be the set of all
$(a_i)_{1\leq i\leq m}\in R(T_{comp})^m$ such that
$a_i\equiv a_j \pmod{(1-e^{\chi})}$ whenever the maximal dimensional
cones $\sigma_i$ and $\sigma_j$ share a wall $\sigma_i\cap \sigma_j$
in $\Delta$ and $\chi\in (\sigma_i\cap \sigma_j)^{\perp}\cap M$. Then
$\mathcal{A}$ is an $R(T_{comp})$-algebra where $R(T_{comp})$ is
identified with the subalgebra of $\mathcal{A}$ consisting of the
diagonal elements $(a,a,\ldots, a)$.  \bth\label{mth1} (see Theorem
\ref{main1}) The ring $K^0_{T_{comp}}(X)$ is isomorphic to
$\mathcal{A}$ as an $R(T_{comp})$-algebra.\eeth Let
\[PLP(\Delta):=\{(f_{\sigma}) \in \prod_{\sigma\in \Delta}
  {\mathbb{Z}[M/M\cap \sigma^{\perp}]}~\mid~f_{\sigma'}\mid_{\sigma}=
  f_{\sigma}~ \mbox{whenever}~ \sigma\preceq \sigma' \in \Delta\}.\]
Then $PLP(\Delta)$ is a ring under pointwise addition and
multiplication and is called the ring of piecewise Laurent polynomial
functions on $\Delta$. The canonical map
$R(T_{comp})=\mathbb{Z}[M]\ra PLP(\Delta)$ which sends $f$ to the
constant tuple $(f)_{\sigma\in \Delta}$ gives $PLP(\Delta)$ the
structure of an $R(T_{comp})$-algebra.

\bth\label{mth2} (see Theorem \ref{main2}) The ring
$K^0_{T_{comp}}(X)$ is isomorphic to $PLP(\Delta)$ as an
$R(T_{comp})$-algebra.
\eeth

Theorem \ref{mth1} and Theorem \ref{mth2} are extensions of the
results of Vezzosi and Vistoli \cite[Theorem 6.2 and Theorem 6.4]{vv}
for algebraic equivariant $K$-ring of smooth toric varieties to the
setting of topological equivariant $K$-ring of any complete
$T$-cellular (possibly singular) toric varieties.

Theorem \ref{mth1} and Theorem \ref{mth2} are also extensions of the
results in \cite[Section 3.1, Section 4]{saruma} to the topological
equivariant $K$-ring of any complete $T$-cellular toric variety (not
necessarily simplicial and projective).

\bth\label{mth3}(see Theorem \ref{basis} and Corollary
\ref{multstrconst}) (1) There exist canonical elements
$f_i\in K^0_{T_{comp}}(X)$ such that
$\displaystyle f_i\mid_{x_i}=\prod_{j=1}^{n-k_i} (1-e^{u_{i_j}})$ and
$f\mid_{x_l}=0$ for all $l>i$. Then $f_i$ form a basis for
$K^0_{T_{comp}}(X)$ as an $R(T_{comp})$-module.

(2) Let $\displaystyle f_i\cdot f_j=\sum_{p=1}^m a^{p}_{i,j}\cdot f_p$
where $a^{p}_{i,j}\in R(T_{comp})$.  The structure constants
$a^{l}_{i,j}$ can be determined by the following closed formula
iteratively, by using descending induction on $l$:
\[ a^l_{i,j}=\frac{\Big[f_i\cdot f_j-\sum_{p=l +1}^{m}a^{p}_{i,j}\cdot
    f_{p}\Big]\mid_{x_l}}{\prod_{r=1}^{n-k_l}(1-e^{-u_{l_r}})}\]for
$1\leq l\leq m$.  \eeth

The above theorem is a generalization of \cite[Theorem 6.9]{gon1} on
the equivariant cohomology of $\mathbb{Q}$-filtrable GKM varieties to
the setting of topological equivariant $K$-ring of $T$-cellular toric
varieties. Also see \cite[Proposition 2.13]{ty} for similar results on
the equivariant cohomology of GKM varieties.

\section{Bialynicki-Birula decomposition and cellular 
  varieties}\label{cellular varieties}

In this section we recall the notions of a filtrable Bialynicki-Birula
cellular decomposition and a $T$-cellular variety (see \cite[Section
3]{Br2} and \cite[Definition 2.1]{gon2}).

Let $X$ be a normal complex algebraic variety with the action of the
complex algebraic torus $T$. We assume that the set of $T$-fixed
points $X^{T}$ is finite. Let $\lambda$ be a generic one-parameter
subgroup of $T$ i.e.  $X^{\lambda}=X^{T}=\{ x_1,\ldots, x_m\}$. Let
\be\label{ps}\displaystyle X_{+}(x_i, \lambda)=\{ x\in X\mid
\lim_{t\ra 0}\lambda(t)x ~ \mbox{exists and is equal to}~x_i\}.\ee
Then $X_{+}(x_i,\lambda)$ is a locally closed $T$-invariant subvariety
of $X$ and is called the Bialynicki-Birula {\it plus} cell.


\bdefe\label{filtrable} The $T$-variety $X$ is called {\em
  filtrable} if it satisfies the following conditions:

(i) $X$ is the union of its plus cells $X_{+}(x_i,\lambda)$ for
$1\leq i\leq m$ with respect to a fixed generic one parameter subgroup
$\lambda$.

(ii) There exists a finite decreasing sequence of $T$-invariant closed
subvarieties of $X$
$X=Z_1\supset Z_2 \cdots\supset  Z_{m} \supset Z_{m+1}= \emptyset $
such that $Z_i\setminus Z_{i+1}=Y_i:=X_{+}(x_i,\lambda)$ for
$1\leq i\leq m$. In particular,
$\displaystyle \overline{Y_i}\subseteq Z_{i}=\bigcup_{j\geq i} Y_j$.
\edefe

\bdefe\label{cellular} Let $X$ be a normal complex algebraic variety
with an action of a complex algebraic torus $T$ such that
$X^T:=\{x_1,\ldots, x_m\}$ is finite. We say that $X$ is {\em $T$-cellular}
if $X$ is filtrable with respect to some fixed generic one-parameter
subgroup $\lambda$ and in addition each plus cell
$Y_i:=X_+(x_i,\lambda)$ is $T$-equivariantly isomorphic to a complex
affine space $\mathbb{C}^{n_i}$ on which $T$-acts linearly for
$1\leq i\leq m$.\edefe

\brem\label{chambers} Note that
$\lambda\in \mbox{Hom}(\mathbb{C}^*,T)=N\simeq \mathbb{Z}^n$.  Let
$\mathfrak{t}:=\mbox{Lie}(T)\simeq N\otimes_{\mathbb{Z}}
\mathbb{C}\simeq \mathbb{R}^{2n}$ as a real vector space. Under the
action of $T$ on $X$ we can consider the stabilizers and their Lie
algebras. A generic one-parameter subgroup $\lambda$ lies in a chamber
which is a connected component of
\[ \mathfrak{t}\setminus \bigcup \{\mbox{Lie algebras of the
    stabilizers of subtori of codimension $\geq 1$}\}.\] Indeed, the
Bialynicki-Birula cell decomposition depends on the chamber in which
$\lambda$ belongs. For example, $-\lambda$ is a generic one-paremeter
subgroup in the opposite chamber and gives rise to the minus
Bialynicki-Birula cells $X_{-}(x_i, \lambda)=X_{+}(x_i,-\lambda)$ by
reversing the direction of the $T$-action (see \cite[Section
3.1]{Br2}). We follow \cite[Section 3.1]{Br2} for the notion of a
generic one-parameter subgroup. Also see \cite[Section 4.1]{car} for
the existence of a generic one-parameter subgroup for a projective
variety, where it is termed as a {\it regular} one-parameter
subgroup. \erem


\brem\label{depends} We note here that the property of being
$T$-cellular also depends on the chamber in which $\lambda$ lies. See
Remark \ref{afterexample1} after Example \ref{example1} where
$T$-cellularity fails for the toric variety $X(\Delta)$ by choosing a
generic one-parameter subgroup lying in a different chamber. However,
by varying $\lambda$ in the same chamber the Bialynicki-Birula cell
decomposition as well as the property of being $T$-cellular does not
change. \erem

\brem\label{rational cellular} More generally, $X$ is called a {\em
  rational} $T$-cellular variety if $X$ is filtrable with respect to a
generic one-parameter subgroup $\lambda$ such that the
Bialynicki-Birula plus cells $Y_i$ for $1\leq i\leq m$ are {\em
  rationally smooth} i.e.
$H^{2n_i}(Y_i, Y_i\setminus \{x_i\};\mathbb{Q})=\mathbb{Q}$ and
$H^m(Y_i, Y_i\setminus \{x_i\};\mathbb{Q})=0$ for $m\neq 2n_i$ where
$n_i=\dim(Y_i)$ (see \cite[Definition 3.4]{gon1} and \cite{Br4}). In
particular, the $\mathbb{Q}$-filtrable varieties considered in
\cite[Defionition 4.6]{gon1} are rational cellular (see \cite[Theorem
4.7]{gon1}). The $T$-cellular varieties defined above in Definition
\ref{cellular} are rational $T$-cellular varieties where the
Bialynicki-Birula cells are all smooth.  \erem




\section{Cellular Toric varieties}\label{cellulartv}

We begin by fixing some notations and conventions.

Let $X=X(\Delta)$ be the toric variety associated to a fan $\Delta$ in
the lattice $N\simeq \bz^n$. Let $M:=\mbox{Hom}(N,\bz)$ be the dual lattice
of characters of $T$. Let $\{v_1,\ldots,v_d\}$ denote the set of
primitive vectors along the edges
$\Delta(1):=\{\rho_1,\ldots,\rho_d\}$. Let $V(\gamma)$ denote the
orbit closure in $X$ of the $T$-orbit $O_{\gamma}$ corresponding to
the cone $\gamma\in\Delta$. Let $S_{\sigma}=\sigma^{\vee}\cap M$ and
$U_{\sigma}:=\mbox{Hom}_{sg}(S_{\sigma},\mathbb{C})$ denote the $T$-stable
open affine subvariety corresponding to a cone $\sigma\in
\Delta$. Here $\mbox{Hom}_{sg}$ denotes semigroup homomorphisms.

We further assume that all the maximal cones in $\Delta$ are
$n$-dimensional, in other words $\Delta$ is {\it pure}. The $T$-fixed
locus in $X$ consists of the set of $T$-fixed points
\be\label{fp}\{x_{1}, x_{2}\ldots,x_{m}\}\ee corresponding to the set
of maximal dimensional cones
\be\label{maxc}\Delta(n):=\{\sigma_1,\sigma_2,\ldots,\sigma_m\}.\ee

Choose a generic one parameter subgroup $\lambda_v\in X_{*}(T)$
corresponding to a {\it generic} $v\in |\Delta|\cap N$ which is
outside the hyperplanes spanned by the $(n-1)$-dimensional cones, so
that (\ref{fp}) is the set of fixed points of $\lambda_v$ (see
\cite[\S3.1]{Br2}).  For each $x_{i}$, we have the plus cell
$Y_i:=X_+(x_{i},\lambda_v)$ (see \eqref{ps}).

Consider a face $\gamma$ of $\sigma_i$ satisfying the property that
the image of $v$ in $N_{\br}/\br\gamma$ is in the relative interior of
$\sigma_i/\br\gamma$. Since the set of such faces is closed under
intersections, we can choose a minimal such face of $\sigma_i$ which
we denote by $\tau_i$. We have \be\label{pu}
Y_i=\bigcup_{\tau_i\subseteq \gamma\subseteq \sigma_i}O_{\gamma}.\ee
Moreover, for a fixed {\it generic} vector $v\in|\Delta|\cap N$,
\[\displaystyle{X=\bigcup_{i=1}^m Y_i}\] where $Y_i$,
$1\leq i\leq m$ are the cells of the {\it Bialynicki-Birula cellular
  decomposition} of the toric variety $X$ corresponding to
the one-parameter subgroup $\lambda_v$ (see \cite[Lemma 2.10]{hs}).

Furthermore, the Bialynicki-Birula decomposition for $X$ corresponding
to $\lambda_v$ is {\it filtrable} if
\be\label{clstar}\bar{Y_{i}}\subseteq \bigcup_{j\geq i} Y_{j}\ee for
every $1\leq i\leq m$.  It can be seen that \eqref{clstar} is
equivalent to the following combinatorial condition in $\Delta$:
\be\label{star} \tau_i\subseteq \sigma_j~~\mbox{implies}~~i\leq j.\ee
Now, from \eqref{pu} it follows that
$Y_i=V(\tau_i) \cap U_{\sigma_i}$. Let $k_i:=\dim(\tau_i)$ for
$1\leq i\leq m$. Thus $Y_i$ is a $T$-stable affine open set in the
toric variety $V(\tau_i)$, corresponding to the maximal dimensional
cone $\bar{\sigma_i}$ in the fan
$(\mbox{star}(\tau_i), N(\tau_i)=N/N_{\tau_i})$ (see \cite[Chapter
3]{f}). Thus $Y_i$ is isomorphic to the complex affine space
$\mathbb{C}^{n-k_i}$ if and only if the $(n-k_i)$-dimensional cone
$\bar{\sigma_i}:=\sigma_i/N_{\tau_i}$ is a smooth cone in
$(\mbox{star}(\tau_i), N(\tau_i))$ for $1\leq i\leq m$.

It follows by \eqref{clstar} that under the above conditions on
$\Delta$, namely \eqref{star} and that $\bar{\sigma_i}$ is a smooth cone in
$\mbox{star}(\tau_i)$, $\displaystyle Z_i:=\bigcup_{j\geq i} Y_j$ are
$T$-stable closed subvarieties, which form a chain
\be\label{strata}X=Z_1\supseteq Z_2\supseteq \cdots \supseteq
Z_m=V(\tau_m)\ee such that
$Z_{i}\setminus Z_{i+1}=Y_i\simeq \mathbb{C}^{n-k_i}$ for
$1\leq i\leq m$.

Since $\bar{\sigma_i}$ is a smooth cone of dimension $n-k_i$ in
$\mbox{star}(\tau_i)$, the primitive vectors
$\bar{v}_{i_1},\ldots, \bar{v}_{i_{n-k_i}}$ in $N(\tau_i)$ along the
edges of $\bar{\sigma_i}$ form a basis for the lattice $N(\tau_i)$.

Now, $\bar{\sigma_i}^{\vee}$ is generated by the dual basis
$u_{i_1},\ldots, u_{i_{n-k_i}}$ for the lattice
$M(\tau_i)=M\cap \tau_i^{\perp}$. Thus any
$x\in Y_i=U_{\bar{\sigma_i}}\simeq \mathbb{C}^{n-k_i}$ can be uniquely
expressed as $(x(u_{i_j}))_{j=1}^{n-k_i}$. Moreover, for $t\in T$ and
$x\in Y_i$, $t\cdot x$ can be expressed as
$(t(u_{i_j})x(u_{i_j}))_{j=1}^{n-k_i}$. Thus $Y_i$ is
$T$-equivariantly isomorphic to the $(n-k_i)$-dimensional complex
$T$-representation
$\displaystyle{\rho_i=\bigoplus_{j=1}^{n-k_i} \rho_{i_j}}$ where
$\rho_{i_j}$ is the $1$-dimensional subrepresentation of $\rho_i$
corresponding to the character $\chi_{i_j}$ where
$\chi_{i_j}(t)=t(u_{i_j})$ for $t\in T$, $1\leq j\leq n-k_i$ and
$1\leq i\leq m$.




We summarize the above discussion in the following theorem which gives
a combinatorial characterization on $\Delta$ for $X(\Delta)$ to be a
$T$-cellular toric variety (see Definition \ref{cellular}).

Fix a generic $v\in |\Delta|\cap N$ such that $\lambda_v$ is a generic
one-parameter subgroup of $T$.

\bth\label{combcell} The toric variety $X(\Delta)$ is $T$-{\it
  cellular} with respect to the Bialynicki-Birula plus cell
decomposition corresponding to $\lambda_v$ if and only if the
following combinatorial conditions hold in $\Delta$:

\noindent
(i) $\Delta$ admits an ordering of the
maximal dimensional cones
\be\label{ordering}\sigma_1<\sigma_2<\cdots<\sigma_m\ee such that the
distinguished faces $\tau_i\subseteq \sigma_i$ defined above satisfy
the following property: \be\label{star1} \tau_i\subseteq
\sigma_j~~\mbox{implies}~~i\leq j 
\ee
\noindent
(ii) $\bar{\sigma_i}:=\sigma_i/N_{\tau_i}$ is a smooth cone in the fan
$(\mbox{star}(\tau_i), N(\tau_i))$ for $1\leq i\leq m$.
\eeth

\bdefe\label{cellular fan} We shall call any fan satisfying the
conditions (i) and (ii) of Theorem \ref{combcell} with respect to $v$
a {\em cellular fan}. \edefe

\brem\label{tau_1} We remark here that we have assumed implicitly that
$v$ belongs to the relative interior of the maximal dimensional cone
$\sigma_1$ with respect to the ordering \eqref{ordering}, so that
$\tau_1=\{0\}$. In particular, if $X(\Delta)$ is cellular with respect
to the one-parameter subgroup $\lambda_v$, then $\sigma_1$ is a smooth
cone in $\Delta=\mbox{star}(\tau_1)$, so that
$Y_1=U_{\sigma_1}\cap V(\tau_1)\simeq \mathbb{C}^n$, the open
Bialynicki-Birula cell corresponding to the fixed point $x_1$, is
smooth in $X=V(\tau_1)$.  \erem

%


\brem When $X=X(\Delta)$ is a smooth projective $T$-toric variety,
then the corresponding fan $\Delta$ is smooth and polytopal. Recall
that we have a filtrable Bialynicki-Birula cellular decomposition of
$X$ corresponding to a generic one-parameter subgroup $\lambda_v$ on
$T$ (see \cite{BB1} and \cite{Br2}).  This in turn gives a generic
height function on the polytope $P$ dual to $\Delta$, which comes from
the linear form $\langle ~ ,v\rangle$ on $M$ defined by $v\in N$. Let
$u(\sigma_i)$ denote the vertex of the polytope corresponding to the
maximal dimensional cone $\sigma_i$. Then the ordering of the maximal
dimensional cones of $\Delta$ such that
$\langle u_{\sigma_1},v\rangle<\langle u_{\sigma_2},v\rangle<\cdots<
\langle u_{\sigma_m},v\rangle$ satisfies the condition (\ref{star})
(see \cite[page 101]{f}).  Moreover, since $x_{\sigma_i}$ is a smooth
point in $X$ it is smooth in $V(\tau_i)$ for $1\leq i\leq m$,
verifying the cellular condition. In particular, a smooth polytopal
fan is a cellular fan. \erem

The toric variety corresponding to the following fan is a $2$-dimensional
non-smooth, non-complete cellular toric variety.

\beg\label{example1} Let $\Delta$ in
$N=\mathbb{Z}e_1\oplus \mathbb{Z}e_2$ consist of the $2$-dimensional
cones $\sigma_1=\langle e_1, 4e_1+e_2\rangle$,
$\sigma_2=\langle 2e_1+e_2, 4e_1+e_2\rangle$ and
$\sigma_3=\langle e_2, 2e_1+e_2\rangle$ and their faces. We denote the
edge generated by $e_1$ by $\rho_1$, the edge generated by $4e_1+e_2$
by $\rho_2$, the edge generated by $2e_1+e_2$ by $\rho_3$, the edge
generated by $e_2$ by $\rho_4$.  Let $v=5e_1+e_2$. Then $v$ belongs to
the relative interior of $\sigma_1$ so that $\tau_1=\{0\}$.  We now
check that $\tau_2=\langle 2e_1+e_2 \rangle=\rho_3$ and
$\tau_3=\langle e_2\rangle=\rho_4$. Note that $\bar{v}=\bar{3e_1}$ in
$N_{\mathbb{R}}/\mathbb{R}\cdot {(2e_1+e_2)}$. On the other hand
$\bar{\sigma_2}=\sigma_2/\mathbb{R}\cdot {(2e_1+e_2)}$ is a
$1$-dimensional cone generated by the class of
$\bar{4e_1+e_2}=\bar{2e_1}$. Since $\bar{e_1}$ is the primitive vector
along $\bar{\sigma_2}$ it follows that $\bar{v}$ is in the relative
interior of $\bar{\sigma_2}$ in
$\mbox{star}(\langle 2e_1+e_2 \rangle)$.  Moreover, since $v$ is not
in the relative interior of $\sigma_2$ it follows that
$\tau_2=\langle 2e_1+e_2\rangle$.  Similarly the image of
$\bar{v}=\bar{5e_1}$ in $N_{\mathbb{R}}/\mathbb{R} \cdot {e_2}$. Now,
$\bar{\sigma_3}=\sigma_3/\mathbb{R} e_2$ is the $1$-dimensional cone
spanned by $\bar{2e_1+e_2}$ and hence generated by the primitive
vector $\bar{e_1}$.  Thus $\bar{v}$ is in the relative interior of
$\bar{\sigma_3}$ in $\mbox{star}(\langle e_2\rangle)$. Since $v$ is
not in the relative interior of $\sigma_3$ it follows that
$\tau_3=\langle e_2\rangle$. We further observe that the ordering
$\sigma_1<\sigma_2<\sigma_3$ satisfies (\ref{star}). It therefore
follows that $X(\Delta)$ is a $T$-cellular toric variety.\eeg

The below picture  illustrates the above example.

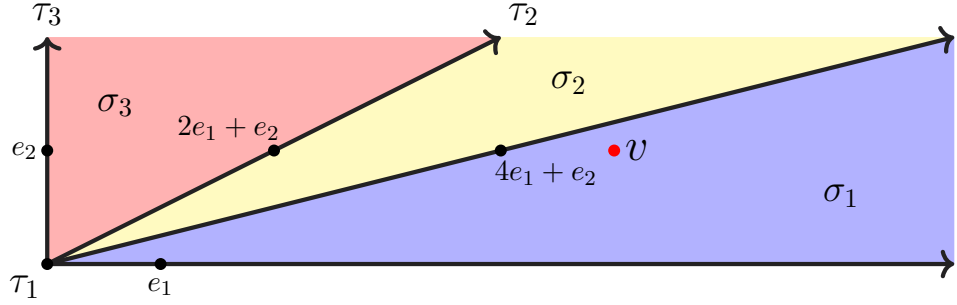
\begin{tikzpicture}[scale=3]
	
	\coordinate (B1) at (0,0);
	\coordinate (B2) at (4,0);
	\coordinate (B3) at (0,1);
	\coordinate (B4) at (4,1);
	\coordinate (B5) at (2,1);
	\coordinate (B6) at (5,1);
	\coordinate (B7) at (-1,-1);
	\coordinate (B8) at (1,-1);

	\definecolor{c1}{RGB}{0,129,188}
	\definecolor{c2}{RGB}{252,177,49}
	\definecolor{c3}{RGB}{35,34,35}

	\fill [blue!30] (B1) -- (B2) -- (B4)  -- cycle;
	\fill [yellow!30] (B1) -- (B4) -- (B5) -- cycle;
	\fill [red!30] (B1) -- (B5) -- (B3) --cycle;

	\draw[c3, ultra thick,->] (B1) -- (B2);
	\draw[c3, ultra thick,->] (B1) -- (B3);
	\draw[c3, ultra thick,->] (B1) -- (B4);
	\draw[c3, ultra thick,->] (B1) -- (B5);

	\fill  (0,0) circle  (.75pt);
	\fill  (.5,0) circle (.75pt);
	\fill  (0,.5) circle (.75pt);
	\fill  (2,.5) circle (.75pt);
	\fill  (1,.5) circle (.75pt);
	\fill[red] (2.5,.5) circle (.75pt);

	\node at (3.5,.3) {\large \(\sigma_1 \)};
	\node at (2.3,.8) {\large  \(\sigma_2 \)};
	\node at (.3,.7) {\large \(\sigma_3 \)};
	\node at (2.6,.5) {\Large \(v \)};
	\node at (-0.1,-0.1) {\large \(\tau_1\)};
	\node at (2.1,1.1) {\large \(\tau_2\)};
	\node at (0,1.1) {\large \(\tau_3\)};
	\node at (.5,-0.1){\small \(e_1\)};
	\node at (-0.1,0.5){\small\(e_2\)};
	\node at (2.2,0.4){\small\(4e_1+e_2\)};
	\node at (0.8,0.6){\small\(2e_1+e_2\)};

	\end{tikzpicture}

        \brem\label{afterexample1} In the above example if we choose
        $v=3e_1+e_2$ in the relative interior of the cone
        $\sigma_2$. Then it can be seen that $\tau_2=\{0\}$. Since
        $\sigma_2$ is not smooth in $\mbox{star}(\tau_2)$ we see that
        $X(\Delta)$ is not $T$-cellular with respect to the one
        parameter subgroup $\lambda_v$. Indeed $\lambda_v$ gives a
        filtrable Bialynicki-Birula cell decomposition of $X(\Delta)$,
        where the Bialynicki-Birula cell corresponding to the
        $T$-fixed point $x_{2}$ is not smooth. Also see Remark
        \ref{tau_1}. \erem

        The following is an example of a complete non-simplicial
        cellular fan in dimension $3$.

        \beg\label{non-simplicial cellular fan} Let $N=\mathbb{Z}^3$
        and $\Delta$ be a fan in $N$ consisting of the following $3$
        dimensional cones and their faces:
        $\sigma_1=\langle e_1, e_2, e_1+e_2+e_3\rangle$,
        $\sigma_2=\langle e_2+e_3, e_1+e_2+e_3, e_2\rangle$,
        $\sigma_3=\langle e_1+e_3, e_1+e_2+e_3, e_1\rangle$,
        $\sigma_4=\langle e_3, e_2+e_3, e_1+e_3, e_1+e_2+e_3\rangle$,
        $\sigma_5=\langle e_1, e_2, -e_1-e_2-e_3\rangle$,
        $\sigma_6=\langle e_1, e_1+e_3, -e_1-e_2-e_3\rangle$,
        $\sigma_7=\langle e_3, e_1+e_3, -e_1-e_2-e_3\rangle$,
        $\sigma_8=\langle e_2+e_3, e_2, -e_1-e_2-e_3\rangle$,
        $\sigma_9=\langle e_3, e_2+e_3, -e_1-e_2-e_3\rangle$. Note
        that $\sigma_i$ for $1\leq i\leq 9$ such that $i\neq 4$, are
        smooth and $\sigma_4$ is not simplicial and hence not smooth.
        Now, by choosing $v=4e_1+3e_2+e_3$ we get that $\tau_1=\{0\}$,
        $\tau_2=\langle e_2+e_3\rangle$,
        $\tau_3=\langle e_1+e_3\rangle$, $\tau_4=\langle e_3\rangle$,
        $\tau_5=\langle -e_1-e_2-e_3\rangle$,
        $\tau_6=\langle e_1+e_3, -e_1-e_2-e_3\rangle$,
        $\tau_7=\langle e_3, -e_1-e_2-e_3\rangle$,
        $\tau_8=\langle e_2+e_3, -e_1-e_2-e_3\rangle$,
        $\tau_9=\sigma_9$. Further, for the ordering
        $\sigma_1<\sigma_2<\cdots <\sigma_8<\sigma_9$ the property
        \eqref{star} holds. Furthermore, since $\sigma_i$ is smooth in
        $\Delta$ for all $1\leq i\leq 9$ such that $i\neq 4$,
        $\bar{\sigma_i}$ is smooth in $\mbox{star}(\tau_i)$. Moreover,
        since $\bar{\sigma_4}=\langle \bar{e_1}, \bar{e_2}\rangle$ in
        $N_{\mathbb{R}}/\mathbb{R}\cdot e_3$, it follows that
        $\bar{\sigma_4}$ is smooth in $\mbox{star}(\tau_4)$. Hence
        $\Delta$ is a cellular fan with respect to $v=4e_1+3e_2+e_3$.
        \eeg

        Henceforth we assume that $X:=X(\Delta)$ is a $T$-cellular
        toric variety with respect to the generic one-parameter
        subgroup $\lambda_v$ of $T$. In particular, we assume that the
        fan $\Delta$ is cellular with respect to $v$.

        \section{Generalized equivariant cohomology theory of
          $T$-cellular varieties }\label{prelimkth}

Let $T\simeq (\mathbb{C}^*)^n$ and
$T_{comp}\simeq (S^1)^n\subseteq T$.




Let $\mbox{h}^*_{G_{comp}}$ denote a generalized $G_{comp}$-equivariant
cohomology theory such that $\mbox{h}^{q}_{G_{comp}}(pt)=0$ for $q$-odd.

\bth\label{k-thcellular} Let $X$ be a $T$-cellular variety (see
Definition \ref{cellular}). Then $\mbox{h}_{T_{comp},c}^0(X)$ is a
free $\mbox{h}^0_{T_{comp}}(pt)$-module of rank $m$ which is the
number of cells. Furthermore, we have
$\mbox{h}_{T_{comp},c}^{-1}(X)=0$.  \eeth

\begin{proof}
  By \cite[Section 1]{mat1} we have the following long exact sequence
  of $T_{comp}$-equivariant generalized cohomology groups which is
  infinite in both directions: \begin{align}\label{les} \cdots\ra
    {\mbox{h}}^{q}_{T_{comp}}(Z_i,Z_{i+1})\ra {\mbox{h}}_{T_{comp},c}^{q}(Z_i)&\ra
                                                                \mbox{h}^{q}_{T_{comp},c}(Z_{i+1})\ra \nonumber\\
                                                              &\ra
                                                                \mbox{h}^{q+1}_{T_{comp}}(Z_i,Z_{i+1})\ra
                                                                \cdots \end{align}
                                                                for
                                                                $1\leq
                                                                i\leq
                                                                m$ and
                                                                $q\in
                                                                \mathbb{Z}$.

Moreover, \cite[Lemma 1.1]{mat1} we have
\begin{align*}{\mbox{h}}^q_{T_{comp}}(Z_i,Z_{i+1})&={\widetilde{\mbox{h}}}^q_{T_{comp}}(Z_i/Z_{i+1})\\ &=\mbox{h}^q_{T_{comp},c }(Y_i)\\
                                                  &=\mbox{h}^{q}_{T_{comp}}(Y_i,
                                                    Y_i\setminus
                                                    \{x_{i}\})= \mbox
                                                    {h}_{T_{comp}}^{q-2n_i}(x_i) \end{align*}
                                                    where
                                                    $Y_i\cong
                                                    \mathbb{C}^{n_i}$
                                                    for
                                                    $1\leq i\leq m$
                                                    . Thus when $q$ is
                                                    odd
                                                    $\mbox{h}^{q}_{T_{comp}}(Z_i,Z_{i+1})=0$.

                                           Moreover, since $Z_m=Y_m$
                                           and $Z_{m+1}=\emptyset$ we
                                           have
                                           $\mbox{h}^0_{T_{comp},c}(Z_m)=
                                           \mbox{h}^0_{T_{comp}}(pt)$
                                           and
                                           $\mbox{h}^{-1}_{T_{comp},c}(Z_m)=\mbox{h}^{-1}_{T_{comp}}(pt)=0$.

    Now, by decreasing induction on $i$ suppose that
    $\mbox{h}_{T_{comp},c}^0(Z_{i+1})$ is a free $\mbox{h}^0_{T_{comp}}(pt)$-module
    for $1\leq i\leq m$ of rank $(m-i)$ and
    $\mbox{h}_{T_{comp},c}^{-1}(Z_{i+1})=0$.  We can start the induction since
    $\mbox{h}^0_{T_{comp},c}(Z_m)=\mbox{h}_{T_{comp}}^0(pt)$ and
    $\mbox{h}^{-1}_{T_{comp},c}(Z_m)=0$.

    Thus from \eqref{les} we get the following split short exact
    sequence of $\mbox{h}_{T_{comp}}^0(pt)$-modules \be \label{eq1} 0\ra
    \mbox{h}_{T_{comp}}^0(Z_i, Z_{i+1})\ra \mbox{h}_{T_{comp},c}^0(Z_i) {\ra} \mbox{h}_{T_{comp},c}^0( Z_{i+1})\ra 0 \ee for
    $1\leq i\leq m$.

    Hence it follows that $\mbox{h}_{T_{comp},c}^0(Z_i)$ is a free
    $\mbox{h}^0_{T_{comp}}(pt)$-module of rank $m-i+1$. Since $X=Z_1$,
    it follows by descending induction on $i$ that
    $\mbox{h}_{T_{comp},c}^0(Z_1)=\mbox{h}_{T_{comp},c}^0(X)$ is a
    free $\mbox{h}^0_{T_{comp}}(pt)$-module of rank $m$.

    Since $\mbox{h}^{-1}_{T_{comp}}(Z_i,Z_{i+1})=0$ as shown above and
    by induction assumption $\mbox{h}^{-1}_{T_{comp},c}(Z_{i+1})=0$,
    it follows from \eqref{les} that
    $\mbox{h}^{-1}_{T_{comp},c}(Z_i)=0$. It follows by descending
    induction on $i$ that
    $\mbox{h}_{T_{comp},c}^{-1}(Z_1)=\mbox{h}_{T_{comp},c}^{-1}(X)=0$.

\end{proof}

\subsection{Equivariant $K$-theory of $T$-cellular varieties}

Let $\mbox{h}^*_{G_{comp}}=K^*_{G_{comp}}$. By Bott periodicity we
have $K_{G_{comp}}^{q}(pt)=K_{G_{comp}}^{-1}(pt)$ for $q$ odd and
$K_{G_{comp}}^q(pt)=K_{G_{comp}}^0(pt)$ for $q$ even. Furthermore,
$K_{G_{comp}}^{-1}(pt)=R(G_{comp})\otimes K^{-1}(pt)=0$ and
$K_{G_{comp}}^0(pt)=R(G_{comp})$. We therefore have the following
result.

\bth\label{k-thcellular} Let $X$ be a $T$-cellular variety (see
Definition \ref{filtrable}). Then $K_{T_{comp},c}^0(X)$ is a free
$R(T_{comp})$-module of rank $m$-which is the number of
cells. Furthermore, we have $K_{T_{comp},c}^{-1}(X)=0$.
\eeth

We have the following corollary for any $T$-cellular variety $X$.

Let $\iota:X^{T_{comp}}=X^{T}\hra X$ denote the inclusion of the set of
$T$-fixed points in $X$.

\bcor\label{localization} The canonical restriction
map
\[ K^0_{T_{comp}}(X)\stackrel{\iota^*}{\lra} K^0_{T_{comp}}(X^{T_{comp}})\cong R(T_{comp})^m\] is
injective where $m=|X^{T_{comp}}|$.  \ecor
\begin{proof}
  Since the prime ideal $(0)$ of $R(T_{comp})$ has support $T_{comp}$, by
  localizing at $(0)$ (see \cite[Proposition 4.1]{segal}), we have
  that
  \[K^0_{T_{comp}}(X)\otimes_{R(T_{comp})} Q(T_{comp})\lra
    K^0_{T_{comp}}(X^{T_{comp}})\otimes_{R(T_{comp})}Q(T_{comp})\] is
  an isomorphism where $Q(T_{comp}):=R(T_{comp})_{\{0\}}$ is the
  quotient field of the integral domain $R(T_{comp})$. This further
  implies that the restriction map
  $\displaystyle K^0_{T_{comp}}(X)\lra K^0_{T_{comp}}(X^{T_{comp}})$ is
  injective, since by Theorem \ref{k-thcellular}, $K^0_{T_{comp}}(X)$
  is a free $R(T_{comp})$-module of rank $m$.\end{proof}


\brem\label{weaker-notion of $T$-cellular} Note that the results of
the Section \ref{prelimkth} hold for a weaker notion of a $T$-cellular
variety, where we can drop the assumption that the $T$ action on the
affine space $\mathbb{C}^{n_i}$ is linear, in Definition
\ref{cellular}. \erem

\section{GKM theory of $X(\Delta)$}\label{GKMtoricsection}

In this section we shall assume that $\Delta$ is a complete cellular
fan (see Definition \ref{cellular fan}).  We shall follow the
notations of Section \ref{cellulartv}.

In this section we shall give a GKM type description for
$K^0_{T_{comp}}(X(\Delta))$ as a
$K^0_{T_{comp}}(pt)=R(T_{comp})$-algebra for a complete $T$-cellular
toric variety $X(\Delta)$.

Since $\mbox{Hom}_{alg}(T,\mathbb{C}^*)=\mbox{Hom}(T_{comp}, S^1)=M$
we have an isomorphism \be\label{iso3'} R(T)\cong R(T_{comp})=
\mathbb{Z}[M]\ee where $\mathbb{Z}[M]=\mathbb{Z}[e^{u} : u\in M]$.  The
isomorphism \eqref{iso3'} identifies the ring
$R(T)=\mathbb{Z}[\mbox{Hom}_{alg}(T,\mathbb{C}^*)]$ of algebraic
representations of $T$ with the ring
$R(T_{comp})=\mathbb{Z}[\mbox{Hom}(T_{comp}, S^1)]$ of topological
representations of $T_{comp}$. More generally, $R(G)=R(G_{comp})$ for
any connected complex reductive algebraic group $G$ and a maximal
compact subgroup $G_{comp}\subseteq G$ (see \cite[p. 272-273]{cg}).

Let $\mathcal{A}$ be the set of all
$(a_i)_{1\leq i\leq m}\in R(T_{comp})^m$ such that
$a_i\equiv a_j \pmod{(1-e^{\chi})}$ whenever the maximal dimensional
cones $\sigma_i$ and $\sigma_j$ share a wall $\sigma_i\cap \sigma_j$
in $\Delta$ and $\chi\in (\sigma_i\cap \sigma_j)^{\perp}\cap M$. In
other words, the $T$-fixed points $x_i$ and $x_j$ lie in the closed
$T$-stable irreducible curve $C_{ij}:=V(\sigma_i\cap \sigma_j)$ in
$X(\Delta)$ and $T$ (and hence $T_{comp}$) acts on $C_{ij}$ through
the character $\chi$. Moreover, $\mathcal{A}$ is a
$R(T_{comp})$-algebra where $R(T_{comp})$ is identified with the
subalgebra of $\mathcal{A}$ consisting of the diagonal elements
$(a,a,\ldots, a)$.

Let $\iota:X^{T}\hra X$ be the inclusion of the set of $T$-fixed
points in $X$.

The following lemma is a refinement of Corollary \ref{localization}
for a complete $T$-cellular toric variety $X(\Delta)$. For a similar
result when $X$ is a full flag variety see \cite[Lemma 2.3]{ml}. Also
see \cite[Chapter 7]{af} for a proof of the following
result in the setting of equivariant cohomology.

\blem\label{localgkm} Let $X=X(\Delta)$ be a complete $T$-cellular
toric variety. Then we have the following inclusion of
$R(T_{comp})$-algebras
\[\iota^*(K^0_{T_{comp}}(X))\subseteq \mathcal{A}.\] \elem
\begin{proof} Let the $T$-fixed points $x_i$ and $x_j$ lie in the
  closed $T$-stable irreducible curve $C_{ij}$ and $T$ (and hence
  $T_{comp}$) act on $C_{ij}$ through the character $\chi$. The
  composition
  \[K^0_{T_{comp}}(X)\stackrel{\iota^*}{\ra} K^0_{T_{comp}}(X^{T})\ra K^0_{T_{comp}}(\{x_i,x_j\})\] equals
  \[K^0_{T_{comp}}(X)\ra
    K^0_{T_{comp}}(C_{ij})\stackrel{\iota_{C_{ij}}^*}{\ra}
    K^0_{T_{comp}}(\{x_i,x_j\}).\]
  Since $X$ is a $T$-toric variety, $C_{ij}=V(\sigma_i\cap \sigma_j)$
  is a $1$-dimensional toric subvariety of $X$ containing two distinct
  $T$-fixed points $x_i$ and $x_j$. In particular, $C_{ij}$ is smooth
  and isomorphic to $\mathbb{P}^1$. We now claim that the image of
  $\iota_{C_{ij}}^*$ consists of pairs of elements
  $(f,g)\in R(T_{comp})\oplus R(T_{comp})$ such that
  $f\equiv g\pmod{(1-e^{-\chi})}$ where $T$ acts on $C_{ij}$ through
  the weight $\chi$.

  The proof of this fact follows by the same arguments as can be found
  in \cite[proof of Lemma 2.3, p. 322]{ml} and \cite[Theorem 1.3,
  p. 274]{u1} so we do not repeat them (also see \cite[Section
  3.4]{Br2}).

  Thus if $x\in K^0_{T_{comp}}(X)$ then $\iota^*(x)_{i}-\iota^*(x)_j$ is
  divisible by $(1-e^{-\chi})$. hence the lemma.
\end{proof}

\bth\label{main1} Let $X=X(\Delta)$ be a complete $T$-cellular toric
variety.  The ring $K^0_{T_{comp}}(X)$ is isomorphic to $\mathcal{A}$
as an $R(T_{comp})$-algebra.

\eeth \begin{proof} By Corollary \ref{localization} and Lemma
  \ref{localgkm} we know that the image of
  $\iota^*(K^0_{T_{comp}}(X))\subseteq \mathcal{A}$. It remains to
  show that $\iota^*(K^0_{T_{comp}}(X))\supseteq \mathcal{A}$. We
  shall follow the methods in \cite[Section 2.5]{ml}.

  Recall from Section \ref{cellulartv} that $k_i=\dim(\tau_i)$ for
  $1\leq i\leq m$. Further, recall that $Y_i$ is $T$-equivariantly
  isomorphic to a $T$-representation $\rho_i=(V_i,\pi_i,x_i)$ where
  $V_i=\mathbb{C}^{n-k_i}$, $\pi_i:V_i\ra x_i$ is the projection and
  $\displaystyle{\rho_i=\bigoplus_{j=1}^{n-k_i} \rho_{i_j}}$ where
  $\rho_{i_j}$ is the one-dimensional subrepresentation of $\rho_i$
  corresponding to the character $\chi_{i_j}$ where
  $\chi_{i_j}(t)=t(u_{i_j})$ for $t\in T$, $1\leq j\leq k_i$ and
  $1\leq i\leq m$.




  Recall that $\bar{\sigma_i}$ is a smooth cone of dimension $n-k_i$
  in $\mbox{star}(\tau_i)$, so that the primitive vectors
  $\bar{v}_{i_1},\ldots, \bar{v}_{i_{n-k_i}}$ in $N(\tau_i)$ along the
  edges of $\bar{\sigma_i}$ form a basis for the lattice $N(\tau_i)$.

  Since $\Delta$ is complete, $\mbox{star}(\tau_i)$ is a complete
  fan. Thus $\bar{\sigma_i}$ shares a wall with $n-k_i$ cones in
  $\mbox{star}(\tau_i)(n-k_i)$, each of which correspond to a maximal
  dimensional cone in $\Delta$ containing $\tau_i$.

Furthermore, we shall assume without loss of generality that for
$1\leq j\leq n-k_i$, the wall
$\langle \bar{v}_{i_1},\ldots, \widehat{\bar{v}_{i_j}},\ldots,
\bar{v}_{i_{n-k_i}} \rangle$ of $\bar{\sigma_i}$, which is orthogonal
to $u_{i_j}$ is $\bar{\sigma_i}\cap \bar{\sigma_{i_j}}$. Thus it
follows that the maximal dimensional cones in $\mbox{star}(\tau_i)$
which share a wall with $\bar{\sigma_i}$ are precisely
$\bar{\sigma_{i_j}}$ for $1\leq j\leq n-k_i$. This is equivalent to
the condition that in $X(\Delta)$ there is a $T$-invariant curve
$C_{ii_j}$ joining the $T$-fixed points $x_i$ and $x_{i_j}$, on which
$T$-acts via the character $u_{i_j}$, for each $1\leq j\leq n-k_i$.


Furthermore, since $\bar{v}_{i_1},\cdots, \bar{v}_{i_{n-k_i}}$ form a
$\mathbb{Z}$-basis for $N(\tau_i)$, the dual basis vectors
$u_{i_1},\ldots, u_{i_{n-k_i}}$ of $M(\tau_i)$ are primitive and
pairwise linearly independent. Since $u_{i_j}$ is primitive and
non-zero for every $1\leq j\leq n-k_i$, the $K$-theoretic equivariant
Euler class $e^{T}(\rho_{i_j})=1-e^{-u_{i_j}}$ is a non-zero divisor
in the unique factorization domain
$R(T_{comp}) =\mathbb{Z}[e^u : u\in M]$. Moreover, since $u_{i_j}$ and
$u_{i_{j'}}$ are linearly independent, it follows that
$1-e^{-u_{i_j}}$ and $1-e^{-u_{i_{j'}}}$ are relatively prime in
$R(T_{comp})$.


By the Thom isomorphism $K^0_{T_{comp}}(V_i)$ is a free
$R(T_{comp})$-module on one generator $g_i$ which restricts under the restriction map
$K^0_{T_{comp}}(V_i)\lra K^0_{T_{comp}}(x_i)$ to $\displaystyle\prod_{j=1}^{n-k_i} (1-e^{-u_{i_j}})$.

Recall from \eqref{eq1} that for every $1\leq i\leq m$ we have a split
exact sequence of $R(T_{comp})$-modules \be\label{splitses} 0\ra
K^0_{T_{comp},c}(Y_i)\stackrel{{\mathtt
    j}_i}{\ra}K^0_{T_{comp}}(Z_i)\stackrel{{\mathtt
    k}_i}{\ra}K^0_{T_{comp}}(Z_{i+1})\ra 0\ee where
$K^0_{T_{comp},c}(Y_i)=K^0_{T_{comp}}(Z_i,Z_{i+1})$.  These induce the
following maps induced by inclusions
\[\begin{array}{lllllllll}
    K^0_{T_{comp}}(x_1) & &K^0_{T_{comp}}(x_2) &  & &K^0_{T_{comp}}(x_{m-1}) & &K^0_{T_{comp}}(x_m)\\
    ~~~~\uparrow~~~&  &~~~~\uparrow~~~&  & &~~~~\uparrow~~~& &~~~~\uparrow ~~~\\
    K^0_{T_{comp}}(Z_1=X)& \stackrel{{\mathtt k}_1}{\ra}&K^0_{T_{comp}}(Z_2)\stackrel{{\mathtt k}_2}{\ra}&\cdots&\ra&K^0_{T_{comp}}(Z_{m-1})&\stackrel{{\mathtt k}_{m-1}}{\ra}&K^0_{T_{comp}}(Z_m)\\
    ~~~~\uparrow {\mathtt j}_{1}~~~& &~~~~\uparrow {\mathtt j}_2~~~& & &~~~~\uparrow {\mathtt j}_{m-1}~~~& &~~~~\uparrow {\mathtt j}_{m}~~~\\
    K^0_{T_{comp},c}(Y_1)&  &K^0_{T_{comp},c}(Y_2)&  & & K^0_{T_{comp},c}(Y_{m-1})& &K^0_{T_{comp},c}(Y_m)  
\end{array}\]

In the above diagram ${\mathtt k}_i$ are all surjective and ${\mathtt j}_i$ are all injective maps.

Thus we can choose $f_i\in K^0_{T_{comp}}(X)$ such that
${\mathtt k}_{i-1}\circ\cdots {\mathtt k}_1(f_i)={\mathtt
  j}_i(g_i)$. Thus
$\displaystyle \iota^*(f_i)_{i}=\prod_{j=1}^{n-k_i} (1-e^{-u_{i_j}})$
and $\iota^*(f_i)_{{l}}=0$ if $l>i$ since
${\mathtt k}_i\circ {\mathtt j}_i=0$.

Let $a\in \mathcal{A}$. We prove by induction on $i$ that if $a_{l}=0$
for $l\geq i$ then $a\in \iota^*(K^0_{T_{comp}}(X))$. This will complete the proof of the theorem.

For $i=1$, since $a=0$ the claim is obvious.

Assume by induction that $a_{l}=0$ for $l\geq i+1$. We know by the
above arguments that $\sigma_i\in \Delta(n)$ shares a wall with
$\sigma_{i_j}\in \Delta(n)$ and $i_j>i$ for $1\leq j\leq
n-k_i$. Moreover, the wall $\sigma_i\cap \sigma_{i_j}$ is orthogonal
to the character $u_{i_j}$. Alternately the $T$-fixed points $x_i$ and
$x_{i_j}$ lie on a $T$-invariant curve $C_{i,i_j}$ in $X(\Delta)$ on
which $T$ acts through the character $u_{i_j}$. Since $i_j>i$, by
induction hypothesis $a_{i_j}=0$ for all $1\leq j\leq n-k_i$. Now,
since $a\in \mathcal{A}$ we have $(1-e^{-u_{i_j}})$ divides
$a_i-a_{i_j}=a_i$ for $1\leq j\leq n-k_i$. This implies that
$\displaystyle\prod_{1\leq j\leq n-k_i} (1-e^{-u_{i_j}})$ divides
$a_i$ since $R(T_{comp})$ is a unique factorization domain and the
factors are relatively prime. Let
$\displaystyle a_i=c_i\cdot \prod_{1\leq j\leq n-k_i}(1-e^{-u_{i_j}})$
where $c_i\in R(T_{comp})$. Thus
$a-\iota^*(c_i\cdot f_i)\in \mathcal{A}$ by Lemma \ref{localgkm}. Also
$(a-\iota^*(c_i\cdot f_i))_{i}=0$. Thus we get
$(a-\iota^*(c_i\cdot f_i))_{l}=0$ for $l\geq i$. By induction there
exists $q\in K^0_{T_{comp}}(X)$ such that
$\iota^*(q)=a-\iota^*(c_i\cdot f_i)$. Thus we get
$a=\iota^*(q+c_i\cdot f_i)$. Hence the theorem.
\end{proof}

\brem\label{composeconstant} Note that we can define the constant map
${\mathrm f}_{ij}:x_{i} \rightarrow x_{j}$ between any
two $T_{comp}$-fixed points of $X$ such that $i<j$, which satisfy
${\mathrm f}_{ik}={\mathrm f}_{ij} \circ {\mathrm f}_{jk}$ for
$m\geq k>j>i \geq 1$. Thus we have the pull-back map of equivariant
$K$-theory
${\mathrm f}_{ij}^*: K^0_{T_{comp}}(x_{j})\rightarrow
K^0_{T_{comp}}(x_{i})$ which satisfies
${\mathrm f}_{ik}^*={\mathrm f}_{ij}^*\circ {\mathrm f}_{jk}^*$
for $m\geq k>j>i\geq 1$. Thus we have the pull-back isomorphisms
\[{\mathrm f}_{im}^*:K^0_{T_{comp}}(x_{m})\rightarrow K^0_{T_{comp}}(x_{i})\] for $1\leq i\leq m$.
This gives $\displaystyle{\prod_{i=1}^m K^0_{T_{comp}}(x_{i})}$
a canonical $K^0_{T_{comp}}(x_{m})$-algebra structure via the
inclusion defined by $({\mathrm f}_{im}^*(a))$ for
$a\in K^0_{T_{comp}}(x_{m})$. Also by identifying
$K_{T_{comp}}^0(x_{i})$ with $R(T_{comp})$ for each
$i=1, \ldots, m$ we see that ${\mathrm f}^*_{ij}$ corresponds to the
identity map of $R(T_{comp})$ for every $1\leq i<j\leq m$ and the
inclusion $({\mathrm f}_{im}^*)$ corresponds to the diagonal
embedding of $R(T_{comp})$ in $R(T_{comp})^m$. \erem
 
\brem The results in this section can alternately be obtained by
applying the theorem of Harada Henriques and Holm \cite[Theorem
3.1]{haheho} as done in \cite{saruma} for the case of cellular toric
orbifolds, but we give a self-contained proof here using methods similar to
those used in \cite{ml}.  \erem

\subsection{Piecewise Laurent polynomial functions on $\Delta$}\label{pl}

Let $X(\Delta)$ be a complete $T$-cellular toric variety.


We first prove the main result of this section with the following
additional combinatorial assumption on $\Delta$. We shall then show
that this assumption holds in particular for any complete fan (see
Lemma \ref{assumcomplete}).

\begin{assum}\label{1-skeleton} The fan $\Delta$ satisfies the
  following combinatorial property. Let $\tau\in \Delta(k)$ and
  $\sigma$ and $\sigma'\in \Delta(n)$
  such that $\tau\preceq \sigma, \sigma'$. There is a sequence
  ${\sigma}={\sigma_1},\ldots, {\sigma_r}={\sigma'}$ of cones in
  $\Delta(n)$ each containing $\tau$, such that
  ${\sigma_{j}}\cap {\sigma_{j+1}}\in \Delta(n-1)$ for
  $1\leq j\leq r-1$.\end{assum}

In this section we show that $K^0_{T_{comp}}(X)$ is isomorphic as an
$R(T_{comp})$-algebra to the ring of piecewise Laurent polynomial
functions on $\Delta$. We follow methods similar to those used in
\cite[Theorem 5.5]{hhrw} and \cite[Theorem 4.2]{saruma}.


Assumption \ref{1-skeleton} is used in the proof of Theorem
\ref{main2}. The geometric interpretation of this assumption is that
for every $\tau\in \Delta$, any two $T$-fixed points in $V(\tau)$, can
be connected by a finite chain of $T$-invariant curves, each of which
is a projective line joining a pair of $T$-fixed points lying in
$V(\tau)$.

Let $\sigma\in \Delta$.  Then $\sigma^{\perp}\cap M$ is a sublattice of
$M$. Let $u_1,\ldots, u_{r}$ be a basis of $\sigma^{\perp}\cap
M$. Let
$${\mathcal K}_{\sigma}:=R(T_{comp})/J_{\sigma}$$
where $J_{\sigma}$ is the ideal generated by
$e^{T_{comp}}({u_{1}}), \ldots ,e^{T_{comp}}({u_{r}})$ where
\[e^{T_{comp}}(u):=1-e^{-u}\] denotes the $T_{comp}$-equivariant
$K$-theoretic Euler class of the $1$-dimensional complex
$T_{comp}$-representation corresponding to $u$. Moreover, since
$u_i$'s are pairwise linearly independent, $e^{T_{comp}}(u_i)$ for
$1\leq i\leq r$ are pairwise relatively prime in $R(T_{comp})$. We
further note that $e^{T_{comp}}(u)\in J_{\sigma}$ for every
$u\in \sigma^{\perp}\cap M$.

For $\sigma, \sigma'\in\Delta$, whenever $\sigma$ is a face of
$\sigma'$ we have the inclusion
$\sigma'^{\perp}\subseteq \sigma^{\perp}$. This gives an inclusion
$J_{\sigma'}\subseteq J_{\sigma}$ which induces a ring homomorphism
$\psi_{\sigma, \sigma'}: {\mathcal K}_{\sigma'}\rightarrow {\mathcal K}_{\sigma}$.

We have the following isomorphism \be\label{isom}
R(T_{comp})/J_{\sigma}\simeq K_{T_{comp}}^0(T/T_{\sigma})=R(T_{\sigma})\ee
(see \cite[Example 4.12]{hhrw}) which further gives the
identifications \be\label{ident} {\mathcal K}_{\sigma}\simeq
K_{T_{comp}}^0(O_{\sigma})=K_{T_{comp}}^0(T/T_{\sigma})=R(T_{\sigma})\ee for
$\sigma\in \Delta$ where $O_{\sigma}$ is the $T$-orbit in $X$
corresponding to $\sigma$ and $T_{\sigma}\subseteq T$ is the stabilizer
of $O_{\sigma}$.



Further, since $X^*(T_{\sigma})=M/M\cap \sigma^{\perp}$ we can
identify $R(T_{\sigma})\simeq \mathcal{K}_{\sigma}$ with
$\mathbb{Z}[M/M\cap \sigma^{\perp}]$ which is the ring of Laurent
polynomial functions on $\sigma$, for $\sigma\in \Delta$.
Furthermore, whenever $\sigma$ is a face of $\sigma'\in\Delta$, the
homomorphism $\psi_{\sigma,\sigma'}$ can be identified with the
natural homomorphism
$\mathbb{Z}[M/\sigma'^{\perp}\cap M]=R(T_{\sigma'})\ra
\mathbb{Z}[M/\sigma^{\perp}\cap M]=R(T_{\sigma})$, given by
restriction of Laurent polynomial functions on $\sigma'$ to $\sigma$.

Let \be\label{plp}PLP(\Delta):=\{(f_{\sigma}) \in \prod_{\sigma\in
  \Delta} {\mathbb{Z}[M/M\cap
  \sigma^{\perp}]}~\mid~\psi_{\sigma,\sigma'} (f_{\sigma'})=
f_{\sigma}~ \mbox{whenever}~ \sigma\preceq \sigma' \in \Delta.\}\ee
Then $PLP(\Delta)$ is a ring under pointwise addition and
multiplication and is called the ring of piecewise Laurent polynomial
functions on $\Delta$. Moreover, we have a canonical map
$R(T_{comp})=\mathbb{Z}[M]\ra PLP(\Delta)$ which sends $f$ to the
constant tuple $(f)_{\sigma\in \Delta}$. This gives $PLP(\Delta)$ the
structure of $R(T_{comp})$-algebra.

The fan $\Delta$ corresponds to a small category $\mathcal{C}(\Delta)$
with objects $\sigma\in \Delta$ and morphisms given by inclusions
$i_{\tau,\sigma}:\tau\subseteq \sigma$ whose initial object is the
zero cone $\{0\}$. Then $\sigma\mapsto {\mathcal K}_{\sigma}$ defines
a contravariant diagram from this small category to the category of
graded commutative $R(T_{comp})$-algebras $\mathcal{C}(alg)$. The
colimit of this diagram exists in $\mathcal{C}(alg)$ and can be
identified with $PLP(\Delta)$ (see \cite[Section 4]{hhrw} and
\cite[Section 4]{saruma}).

\bth\label{main2} The ring $K^0_{T_{comp}}(X)$ is isomorphic to
$PLP(\Delta)$ as an $R(T_{comp})$-algebra.  \eeth{\bf Proof:} Using
the universal property of colimits we prove this by finding compatible
homomorphisms $h_{\sigma}:\mathcal{A}\ra {\mathcal K}_{\sigma}$ for
every $\sigma\in \Delta$. Let $y=(y_i)\in \mathcal{A}$. On the maximal
dimensional cones of $\Delta$ we define $h_{\sigma_j}(y):=y_j$ for
each $1\leq j\leq m$. On the $(n-1)$-dimensional cones we let
$h_{\gamma}(y):=y_{i} ~\mbox{mod}~ J_{\gamma}$ in
${\mathcal K}_{\gamma}$ for any maximal dimensional cone $\sigma_i$
containing $\gamma$.  If $\gamma=\sigma_i\cap \sigma_j$ for $i\neq j$
then $y_{i}\equiv y_{j} \pmod {e^{T_{comp}}(u_{i_j})}$ where
$u_{i_j}\in \gamma^{\perp}\cap M$ and
$e^{T_{comp}}(u_{i_j})\in J_{\gamma}$. Thus $h_{\gamma}$ is well
defined on the $(n-1)$-dimensional cones. Now, let
$\gamma\in \Delta(n-k)$ for $k\geq 2$.  We can similarly define
$h_{\gamma}(y):=y_i~\mbox{mod}~J_{\gamma}$ where
$\gamma\prec\sigma_i\in \Delta(n)$. To show that $h_{\gamma}$ is well
defined we need to show that if $\gamma\in \sigma_j$ for $j\neq i$
then $y_i\equiv y_j ~\pmod {J_{\gamma}}$. Now,
$\gamma\prec \sigma_i\cap \sigma_j$.  By Assumption \ref{1-skeleton}
we have a sequence of cones
$\sigma_i=\sigma_{i_1},\ldots, \sigma_{i_r}=\sigma_j$ in $\Delta(n)$
all containing $\gamma$ such that
$\mbox{dim}(\sigma_{i_t}\cap \sigma_{i_{t+1}})=n-1$ for
$1\leq t\leq r-1$. Since
$J_{\sigma_{i_t}\cap \sigma_{i_{t+1}}}\subseteq J_{\gamma}$ we have
$y_{i_t}-y_{i_{t+1}}\in {J_{\gamma}}$ for $1\leq t\leq r-1$. Thus
$\displaystyle{y_i-y_j=\sum_{t=1}^{r-1} y_{i_t}-y_{i_{t+1}}\in
  {J_{\gamma}}}$.

Now to check the compatibility of $h_{\gamma}$ for $\gamma\in
\Delta$. For $\gamma\prec \gamma'$ we have
$\psi_{\gamma,\gamma'} : {\mathcal K}_{\gamma'}\rightarrow {\mathcal
  K}_{\gamma}$. We need to verify that
$\psi_{\gamma,\gamma'}\circ h_{\gamma'}=h_{\gamma}$. For
$y\in \mathcal{A}$, $h_{\gamma'}(y)=y_i~\mbox{mod}~J_{\gamma'}$ for
$\gamma'\prec\sigma_i\in \Delta(n)$. Now, by definition of
$\psi_{\gamma,\gamma'}$ we get
$\psi_{\gamma,\gamma'}\circ h_{\gamma'}(y)=y_i~
\mbox{mod}~J_{\gamma}=h_{\gamma}(y)$. This proves the compatibility of
$\{h_{\gamma}\}_{\gamma\in \Delta}$.Thus $\{h_{\gamma}\}$ induces a
well defined ring homomorphism $h$ from $\mathcal{A}$ to
$PLP(\Delta)$. We shall now show that $h$ is an isomorphism. Let
$y\neq y'\in \mathcal{A}$. Then $y_i\neq y'_i$ for at least one
$1\leq i\leq m$. Thus $h_{\sigma_i}(y)\neq h_{\sigma_i}(y')$ in
$R(T_{comp})$. Hence $h$ is injective. Now, let
$(a_{\gamma})_{\gamma\in \Delta}\in PLP(\Delta)$. This in particular
determines an element $(a_{\sigma_i})\in (R(T_{comp}))^m$. We claim
that $(a_{\sigma_i})\in \mathcal{A}$. Let $\sigma_i$ share a wall
$\gamma=\sigma_i\cap \sigma_j$ with $\sigma_j$. Then we see that
$\psi_{\gamma, \sigma_i} (a_{\sigma_i})=a_{\gamma}=\psi_{\gamma,
  \sigma_j}(a_{\sigma_j})$. It implies that
$a_{\sigma_i}-a_{\sigma_j}\equiv 0\pmod {J_{\gamma}}$ whenever
$\sigma_i$ and $\sigma_j$ share a wall $\gamma$ in $\Delta$. Since
$J_{\gamma}$ is generated by $e^{T_{comp}}(u_{i_j})$ where $u_{i_j}$
generates the lattice $\gamma^{\perp}\cap M$, it follows that
$e^{T_{comp}}(u_{i_j})$ divides $a_{\sigma_i}-a_{\sigma_j}$. This
implies by \eqref{main1} that $(a_{\sigma_i})\in \mathcal{A}$ proving
the surjectivity of $h$. $\Box$

\subsection{Strongly connected polytopal complexes}\label{scpc}

We recall the following concepts on polytopal complexes from \cite{z},
\cite{brbe} and \cite{hachimori}. For convenience we borrow the
notation from \cite{brbe}.

A {\em polytopal complex} is a finite, nonempty collection
$\mathcal{P}$ of polytopes (called the faces of $\mathcal{P}$) in some
Euclidean space $\mathbb{R}^k$, such that: (1) if $\sigma$ is in
$\mathcal{P}$ then all the faces of the polytope $\sigma$ are elements
of $\mathcal{P}$; (2) the intersection of any two polytopes $\sigma$
and $\tau$ of $\mathcal{P}$ is a face of both $\sigma$ and $\tau$.


The maximal faces of a polytopal complex under inclusion are called
{\em facets}, and the maximal proper subfaces of the facets under
inclusion are called {\em ridges}.

The dimension of $\mathcal{P}$ is the largest dimension of
a polytope of $\mathcal{P}$.

Pure complexes are polytopal complexes for which all facets have the
same dimension.  The dual graph of a pure $d$-dimensional polytopal
complex $\mathcal{P}$ is a graph whose vertices correspond to the
facets of $\mathcal{P}$ and two vertices of the graph are connected by
an edge if and only if the corresponding facets of $\mathcal{P}$ have
a common ridge.

A pure $d$-dimensional polytopal complex is called {\em strongly
  connected} if its dual graph is connected. In other words, a pure
polytopal complex $\mathcal{P}$ is strongly connected if any two
facets $F$ and $G$ of $\mathcal{P}$ have a sequence
$F=F_1, F_2,\ldots, F_k=G$ of facets such that $F_i$ and $F_{i+1}$
have a common ridge for $1\leq i\leq k-1$. A {\em pseudomanifold} is a
strongly connected polytopal complex for which every ridge is
contained in at most two facets. The underlying space $|\mathcal P|$
of a polytopal complex $\mathcal P$ is the union of all its faces. In
this case, $\mathcal{P}$ is called a polytopal subdivision of
$|\mathcal P|$. We shall call a $d$-dimensional polytopal complex
whose underlying space is homeomorphic to a topological manifold
(compact and connected) as a $d$-manifold. It is known that a
polytopal subdivision of any connected $d$-manifold (with or without
boundary) is pure and strongly connected.



\blem\label{assumcomplete} Any complete fan $\Delta$ in $N$ satisfies
Assumption \ref{1-skeleton}.\elem
\begin{proof} Let $\Delta$ be a complete fan in
  $N\simeq \mathbb{Z}^n$. Then
  $S_{N}:=N_{\mathbb{R}}\setminus \{0\}/\mathbb{R}_{>0}$ is
  homeomorphic to the $(n-1)$-sphere. Consider the projection
  $\uppi:N_{\mathbb{R}}\setminus \{0\}\ra S_{N}$. Then for a strongly
  convex rational polyhedral cone $\sigma\in \Delta$,
  $\uppi(\sigma\setminus \{0\})$ is a spherical polytopal cell of
  dimension $\dim(\sigma)-1$ in $S_N$. Thus a finite complete fan in
  $N$ gives rise to a finite convex spherical polytopal cell
  decomposition of the sphere $S_{N}$ (see \cite[p. 52]{o}). In
  particular, when $\Delta$ is simplicial then
  $\uppi(\sigma\setminus\{0\})$ for $\sigma\in \Delta$ gives a
  triangulation of $S_{N}$ into spherical simplices. Moreover, for any
  complete fan $\Delta$, the polytopal complex $\mathcal{P}_{\Delta}$
  formed by $\uppi(\sigma\setminus \{0\})$ is pure of dimension
  $n-1$. Now, since the underlying space of the polytopal complex
  formed by $\uppi(\sigma\setminus\{0\})$ is an $(n-1)$-sphere, it is a
  pure and strongly connected $(n-1)$-complex. Furthermore, since
  $\Delta$ is complete, for $\tau\in \Delta(k)$,
  $(\mbox{star}(\tau), N(\tau))$ is also complete. In particular,
  $\uppi(\sigma\setminus \tau)$ as $\sigma$ varies over
  $\mbox{star}(\tau)$, form the link of $\uppi(\tau)$. This gives a
  polytopal cell decomposition of an $(n-k-1)$-dimensional sphere and
  is therefore a pure and strongly connected $(n-k-1)$-complex. Note
  that every $n$-dimensional cone of $\Delta$ corresponds to a unique
  facet of $\mathcal{P}_{\Delta}$ and every $(n-1)$-dimensional cone
  of $\Delta$ corresponds to a ridge of $\mathcal{P}_{\Delta}$. Thus
  the strong connectivity of $\mathcal{P}_{\Delta}$ and
  $\mathcal{P}_{\mbox{star}(\tau)}$ for every $\tau\in \Delta$ implies
  that $\Delta$ satisfies Assumption \ref{1-skeleton}.
\end{proof}


\section{Basis for $K^0_{T_{comp}}(X)$ as an $R(T_{comp})-module$}\label{basis1}
Let $X=X(\Delta)$ be a complete $T$-cellular toric variety. In this
section we shall give a basis for $K^0_{T_{comp}}(X)$ as a
$R(T_{comp})$-module.

We shall follow the notations developed in Theorem \ref{main1} and its
proof. In particular, we recall the maps
${\mathtt k}_i:K^0_{T_{comp},c}(Y_i)\ra K^0_{T_{comp}}(Z_i)$ and
${\mathtt j}_i: K^0_{T_{comp}}(Z_i)\ra K^0_{T_{comp}}(Z_{i+1})$ induced by
inclusions and the split exact sequence \eqref{splitses} of
$R(T_{comp})$-modules.

Let $f_i\in K^0_{T_{comp}}(X)$ such that
${\mathtt k}_{i-1}\circ\cdots \circ {\mathtt k}_1(f_i)={\mathtt j}_i(g_i)$. Thus
$\displaystyle \iota^*(f_i)_{i}=\prod_{i=1}^{n-k_i} (1-e^{-u_{i_j}})$ and $\iota^*(f_i)_{{l}}=0$
if $l>i$ since ${\mathtt k}_i\circ {\mathtt j}_i=0$.

We shall let \[f\mid_{x_l}:=\iota^*(f)_{{l}}\] denote the restriction
of an element $f\in K^0_{T_{comp}}(X)$ to the $T$-fixed point $x_l$ for
$1\leq l\leq m$.

\bth\label{basis}  The elements $f_i$ for $1\leq i\leq m$ form a basis for
$K^0_{T_{comp}}(X)$ as a $R(T_{comp})$-module. 

\eeth
\begin{proof} We shall prove by downward induction on $i$ that the set
  $\{\bar{f_i}\}_{r=i}^m$ span $K^0_{T_{comp}}(Z_i)$ as
  $R(T_{comp})$-module. Here
  $\displaystyle\bar{f_r}:={\mathtt k}_{i-1}\circ\cdots \circ {\mathtt k}_1(f_r)\in
  K^0_{T_{comp}}(Z_i)$.  The proposition will follow since $Z_1=X$.

  When $r=m$, $Z_{m}=Y_m$ since $Z_{m+1}=\emptyset$. Thus the
  inclusion ${\mathtt j}_m$ is an isomorphism. Now, $g_m$ spans for
  $K^0_{T_{comp}}(Y_m)$. Hence
  ${\mathtt j}_m(g_m)={\mathtt k}_{m-1}\circ\cdots \circ {\mathtt k}_{1}(f_m)=\bar{f_m}$ spans
  $K^0_{T_{comp}}(Z_m)$.

  Suppose by induction that the claim holds for
  $K^0_{T_{comp}}(Z_{i+1})$. We shall now prove it for
  $K^0_{T_{comp}}(Z_{i})$.

  By induction assumption $\bar{f_r}:={\mathtt k}_{i}\circ\cdots \circ {\mathtt k}_1(f_r)$ for
  $i+1\leq r\leq m$ span $K^0_{T_{comp}} (Z_{i+1})$.

  Let $f\in K^0_{T_{comp}}(Z_i)$. Then
  $\displaystyle {\mathtt k}_i(f)=\sum_{r=i+1}^ma_r\cdot {\mathtt
    k}_{i}\circ\cdots \circ {\mathtt k}_1(f_r)$ where
  $a_r\in R(T_{comp})$ for $i+1\leq r\leq m$.
  
  Thus
  $\displaystyle f-\sum_{r=i+1}^m a_r\cdot {\mathtt
    k}_{i-1}\circ\cdots \circ {\mathtt k}_1(f_r)\in
  \mbox{ker}({\mathtt k}_i)$. From \eqref{splitses} we get that
  $\displaystyle f-\sum_{r=i+1}^m a_r\cdot {\mathtt k}_{i-1}\circ\cdots \circ
  {\mathtt k}_1(f_r)={\mathtt j}_i(a_ig_i)$ for some $a_i\in R(T_{comp})$. Here we recall
  that $g_i$ denotes the Thom class of the $T_{comp}$-equivariant
  vector bundle $V_i$ in $K^0_{T_{comp}}(V_i)$. Thus by definition of
  $f_i$ it follows that
  $\displaystyle f-\sum_{r=i+1}^m a_r\cdot {\mathtt k}_{i-1}\circ\cdots \circ
  {\mathtt k}_1(f_r)=a_i \cdot {\mathtt k}_{i-1}\circ\cdots \circ {\mathtt k}_1(f_i)$. Thus it
  follows that ${\mathtt k}_{i-1}\circ\cdots \circ {\mathtt k}_1(f_r)$ for $i\leq r\leq m$
  span $K^0_{T_{comp}}(Z_i)$. By induction it follows that
  $\{f_i\}_{i=1}^m$ span $K^0_{T_{comp}}(X)$ as a free
  $R(T_{comp})$-module.

  It only remains to show that $\{f_r\}_{r=1}^m$ are linearly
  independent.  Suppose $\displaystyle\sum_{r=1}^m a_r\cdot
  f_r=0$. Then restricting to the fixed point $x_m$ it follows that
  $\displaystyle a_m \cdot \prod_{j=1}^{n-k_{m}} (1-e^{-u_{m_j}})=0$
  since $f_r\mid_{x_m}=0$ for $1\leq r<m$ and
  $\displaystyle f_m\mid_{x_m}=\prod_{j=1}^{n-k_{m}}
  (1-e^{-u_{m_j}})$. Since
  $\displaystyle \prod_{j=1}^{n-k_{m}} (1-e^{-u_{m_j}})$ is a non-zero
  divisor in $R(T_{comp})$ this implies that $a_m=0$. Now, by
  induction we assume that $a_r=0$ for $i+1\leq r\leq m$. We shall now
  show that $a_i=0$. The proof will then follow by induction.

  Now, restricting $\displaystyle\sum_{r=1}^{m} a_i\cdot f_i=0$ to
  $x_{i}$ we get
  $\displaystyle a_i\cdot \prod_{j=1}^{n-k_{i}}
  (1-e^{-u_{i_j}})=0$. This is because by induction assumption and by
  the choice of $f_r$, $a_r=0$ for $i+1\leq r\leq m$,
  $f_r\mid_{x_i}=0$ for $1\leq r\leq i-1$ and
  $\displaystyle f_i\mid_{x_i}=\prod_{j=1}^{n-k_{i}}
  (1-e^{-u_{i_j}})$, which is a non-zero divisor in
  $R(T_{comp})$. This implies that $a_i=0$.\end{proof}


\brem\label{knutsontao} The basis $\{f_i, ~~1\leq i\leq m\}$ can be
viewed as a $K$-theoretic analogue of the notion of Knutson-Tao basis
proposed by Tymoczco in \cite[Definition 2.12, Proposition 2.13]{ty}
for equivariant cohomology of GKM varieties. \erem

\subsection{Structure constants}\label{strconst}

We have the following theorem which gives a closed formula for the
coefficients of an element $f\in K^0_{T}(X)$ with respect to the basis
$\{f_i\}_{i=1}^m$.  \bth\label{coefficients} Let
$f\in K^0_{T_{comp}}(X)$. Let \be\label{expression} \displaystyle
f=\sum_{j=1}^m a_j\cdot f_j\ee where $a_j\in R(T_{comp})$ for
$1\leq j\leq m$.  We have the following closed formula for the
coefficients which is determined iteratively
\be\label{E1}a_{i}=\frac{\Big[f-\sum_{j=i +1}^{m}a_{j}\cdot
  f_{j}\Big]\mid_{x_i}}{\prod_{r=1}^{n-k_i}(1-e^{-u_{i_r}})}\ee for
$1\leq i\leq m$.  \eeth
\begin{proof}
  Recall that $f_j\mid_{x_i}=0$ for $i>j$ and
  $\displaystyle f_j\mid_{x_j}=\prod_{r=1}^{n-k_j}(1-e^{-u_{j_r}})$.

  Thus from \eqref{expression} we get that:
  \[f\mid_{x_m}=\sum_{j=1}^m a_j\cdot f_j\mid_{x_m}=a_m\cdot
    \prod_{r=1}^{n-k_m}(1-e^{-u_{m_r}}). \] Since
  $\displaystyle \prod_{r=1}^{n-k_m}(1-e^{- u_{m_r}})$ is a non-zero
  divisor in $R(T_{comp})$ the theorem follows for $i=m$.  Proceeding
  by descending induction, having determined the coefficients
  $a_m,\ldots, a_{i+1}$, we have the following from
  \eqref{expression}:
  \[ (f-\sum_{j=i+1}^m a_j\cdot f_j)\mid_{x_i}=a_i\cdot
    \prod_{r=1}^{n-k_i}(1-e^{-u_{i_r}}).\]

  Since $\displaystyle \prod_{r=1}^{n-k_i}(1-e^{-u_{i_r}})$ is a
  non-zero divisor in $R(T_{comp})$, the theorem follows for $i$.
  
\end{proof}  

We have the following corollary which gives a formula for the
multiplicative structure constants with respect to the basis
$\{f_i\}_{i=1}^m$.

\bcor\label{multstrconst} Let
$\displaystyle f_i\cdot f_j=\sum_{p=1}^m a^{p}_{i,j}\cdot f_p$ where
$a^{p}_{i,j}\in R(T_{comp})$.  Then the structure constants
$a^{l}_{i,j}$ can be determined by the following closed formula
iteratively, by using descending induction on $l$:
\be\label{closedform} a^l_{i,j}=\frac{\Big[f_i\cdot f_j-\sum_{p=l
    +1}^{m}a^{p}_{i,j}\cdot
  f_{p}\Big]\mid_{x_l}}{\prod_{r=1}^{n-k_l}(1-e^{-u_{l_r}})}\ee for
$1\leq l\leq m$.  \ecor
\begin{proof}
The corollary follows readily by letting $f=f_i\cdot f_j$ in Theorem \ref{coefficients}.
\end{proof}

\brem In Theorem \ref{coefficients} and Corollary \ref{multstrconst},
although a priori the expressions on the right hand side of \eqref{E1}
and \eqref{closedform} respectively belong to $Q(T_{comp})$, as a
consequence of Theorem \ref{basis}, it follows that in fact they
belong to $R(T_{comp})$.  \erem

The following is a picture of a $2$-dimensional non-smooth complete
cellular fan $\Delta$ which contains the fan of Example \ref{example1}
as a subfan.

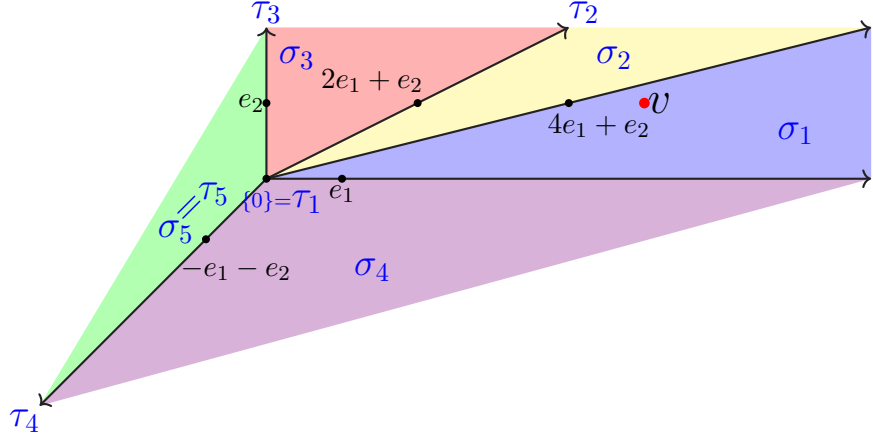
\begin{tikzpicture}[scale=2]
	
	\coordinate (B1) at (0,0);
	\coordinate (B2) at (4,0);
	\coordinate (B3) at (0,1);
	\coordinate (B4) at (4,1);
	\coordinate (B5) at (2,1);
	\coordinate (B6) at (4.5,1);
	\coordinate (B7) at (-1,-1);
	\coordinate (B8) at (1,-1);
        \coordinate (B9) at (-1.5,-1.5);

	\definecolor{c3}{RGB}{35,34,35}

	\fill [blue!30] (B1) -- (B2) -- (B4)  -- cycle;
	\fill [yellow!30] (B1) -- (B4) -- (B5) -- cycle;
	\fill [red!30] (B1) -- (B5) -- (B3) --cycle;
	\fill [green!30] (B1)--(B3)--(B9)--cycle;
        \fill [violet!30] (B1)--(B2)--(B9)--cycle;
	
	\draw[c3,  thick,->] (B1) -- (B2);
	\draw[c3, thick,->] (B1) -- (B3);
	\draw[c3,  thick,->] (B1) -- (B4);
	\draw[c3,  thick,->] (B1) -- (B5);
        \draw[c3,  thick,->] (B1)--(B9);

	\fill  (0,0) circle  (.75pt);
	\fill  (.5,0) circle (.75pt);
	\fill  (0,.5) circle (.75pt);
	\fill  (2,.5) circle (.75pt);
	\fill  (1,.5) circle (.75pt);
        \fill (-.4,-.4) circle (.75pt);
	\fill[red] (2.5,.5) circle (1pt);

	\node at (3.5,.3) {\large \(\textcolor{blue}{\sigma_1} \)};
	\node at (2.3,.8) {\large  \(\textcolor{blue}{\sigma_2} \)};
	\node at (.2,.8) {\large \(\textcolor{blue}{\sigma_3} \)};
        \node at (-.6,-0.35){\large \(\textcolor{blue}{\sigma_5}\)};
        \node at (-.5,-0.2){\large\(\textcolor{blue}{\dequalto}\)};
        \node at (-.35, -0.1){\large \(\textcolor{blue}{\tau_5}\)};
        \node at (0.7,-0.6){\large \(\textcolor{blue}{\sigma_4}\)};
	\node at (2.6,.5) {\Large \(v\)};
	\node at (0.1,-0.15) {\tiny{\textcolor{blue}{\{0\}=}}\large \(\textcolor{blue}{\tau_1}\)};
	\node at (2.1,1.1) {\large \(\textcolor{blue}{\tau_2}\)};
	\node at (0,1.1) {\large \(\textcolor{blue}{\tau_3}\)};
        \node at ( -1.6,-1.6) {\large\(\textcolor{blue}{\tau_4}\)};
	\node at (.5,-0.1){\small\(e_1\)};
	\node at (-0.1,0.5){\small\(e_2\)};
	\node at (2.2,0.35){\small\(4e_1+e_2\)};
	\node at (0.7,0.65){\small\(2e_1+e_2\)};
        \node at (-0.2, -0.6){\small\(-e_1-e_2\)};

	\end{tikzpicture}

        We shall illustrate below the computation of the basis for
        $K^0_{T_{comp}}(X(\Delta))$ for the above $\Delta$.

        \beg\label{Example 3} Now, {$v=5e_1+e_2$}
        belongs to the relative interior of
        {$\sigma_1$} so that
        {$\tau_1=\{0\}$}.  Also,
        {$\tau_2=\langle 2e_1+e_2\rangle$} since
        {$\overline{v}=3\overline{e_1}$} is in the
        relative interior of
        {$\overline{\sigma_2}=\langle \overline{e_1}
          \rangle$} in
        {$N_{\mathbb{R}}/\mathbb{R}\cdot
          {(2e_1+e_2)}$}.  Similarly
        {$\tau_3=\langle e_2 \rangle$} since
        {$\overline{v}=5\overline{ e_1}$} is in the
        relative interior of
        {$\overline{\sigma_3}=\langle
          \overline{e_1}\rangle$} in
        {$N_{\mathbb{R}}/\mathbb{R} \cdot
          {e_2}$}. Furthermore,
        {$\tau_4=\langle -e_1-e_2 \rangle$} since
        {$\bar{v}=4\overline{e_1}$} is in the relative
        interior of
        {$\overline{\sigma_4}=\langle
          \bar{e_1}\rangle$} and {$\tau_5=\sigma_5$}
        since
        {$\overline{v}=5\overline{e_1}=-4\overline{e_2}$}
        is neither in the relative interior of
        {$\overline{\sigma_5}=\langle \overline{e_2}
          \rangle$} in
        {$N_{\mathbb{R}}/\mathbb{R} \langle-e_1-e_2
          \rangle$} nor in the relative interior of
        {$\overline{\sigma_5}=\langle \overline{-e_1}
          \rangle$} in
        $N_{\mathbb{R}}/\mathbb{R} \langle e_2
          \rangle$.

Note that for the {ordering}
  {$\sigma_1<\sigma_2<\sigma_3<\sigma_4<\sigma_5$} the {property
  \eqref{star} holds}.
 Moreover, {$\sigma_i$} is smooth in
  {$\Delta$} for all {$i=1,4,5$}
  corresponding {$\overline{\sigma_i}$} is
  {smooth} in {$\mbox{star}(\tau_i)$}.

  Furthermore, {$\overline{\sigma_2}=\langle \overline{e_1}\rangle$}
  is {smooth} in
  {$(\mbox{star}(\tau_2), N_{\mathbb{R}}/\mathbb{R}\cdot 2e_1+e_2)$}
  and {$\overline{\sigma_3}=\langle \overline{e_1} \rangle$} is
  {smooth} in
  {$(\mbox{star}(\tau_3), N_{\mathbb{R}}/\mathbb{R}\cdot e_2)$}.

  Hence {$\Delta$} is a { cellular fan} with respect to {$v=5e_1+e_2$}.

{$\chi_{12}=e_1^*-4e_2^*$} is {orthogonal to the
  facet} {$\sigma_1\cap \sigma_2$}.

 {$\chi_{14}=e_2^*$} is {orthogonal to the facet} {$\sigma_1\cap \sigma_4$}.

 {$\chi_{23}=e_1^*-2e_2^*$} is {orthogonal to the
  facet} {$\sigma_2\cap \sigma_3$}.

{$\chi_{35}=e_1^*$} is {orthogonal to the facet} {$\sigma_3\cap \sigma_5$}.

{$\chi_{45}=e_1^*-e_2^*$} is {orthogonal to the facet} {$\sigma_4\cap \sigma_5$}.

By Theorem \ref{main1} {$f\in K_{T_{comp}}(X(\Delta))$} {if and only if}\\
  \vspace{0.2cm}
  {$f\mid_{x_1}\equiv f\mid_{x_2} \pmod{1-e^{e_1^*-4e_2^*}}$},\\
   \vspace{0.2cm}
   {$f\mid_{x_1}\equiv f\mid_{x_4} \pmod{1-e^{e_2^*}}$},\\
   \vspace{0.2cm}
   {$f\mid_{x_2}\equiv f\mid_{x_3} \pmod{1-e^{e_1^*-2e_2^*}}$},\\
   \vspace{0.2cm}
   {$f\mid_{x_3}\equiv f\mid_{x_5} \pmod{1-e^{e_1^*}}$},\\
   \vspace{0.2cm}
  {$f\mid_{x_4}\equiv f\mid_{x_5} \pmod{1-e^{e_1^*-e_2^*}}$}.

  By Thorem \ref{basis} we get the following basis for
  $K_{T_{comp}}^0(X(\Delta))$.

 {$f_1=\Big(1-e^{e_2^*}-e^{e_1^*-4e_2^*}+e^{e_1^*-3e_2^*},
  0,0,0,0\Big)$}

{$f_2=\Big(e^{e_1^*-4e_2^*}-e^{e_1^*-2e_2^*} ,
    1-e^{e_1^*-2e_2^*},0,0,0\Big)$}

  {$f_3=\Big(e^{e_1^*-2e_2^*}-e^{e_1^*}
  +1-e^{e_2^*}-e^{e_1^*-4e_2^*}+e^{e_1^*-3e_2^*},
   e^{e_1^*-2e_2^*}-e^{e_1^*}, 1-e^{e_1^*}, 0, 0 \Big)$}

  {$f_4=\Big(1-e^{e_1^*}-e^{e_1^*-2e_2^*}+e^{e_1^*-4e_2^*} ,
    2-e^{e_1^*}-e^{e_1^*-2e_2^*}, 1-e^{e_1^*}, 1-e^{e_1^*-e_2^*},
   0\Big)$}

 {$f_5=\Big(e^{e_1^*}-e^{e_1^*-2e_2^*}+e^{e_1^*-4e_2^*},
                  e^{e_1^*}+1-e^{e_1^*-2e_2^*}, e^{e_1^*}, e^{e_1^*-e_2^*}, 1\Big)$}

\eeg

\brem\label{compareoperationalbasis} In \cite[Example 1.7]{ap} an
example of basis construction for
$\mbox{op}\mathcal{K}_{T}^0(X)=\mbox{Hom}_{R(T)}(\mathcal{K}_0^{T}(X),R(T))$
is discussed for a complete toric variety $X(\Delta)$. More precisely,
for a $T$-cellular toric variety a dual basis
$[\mathcal{O}_{V(\tau_i)}]^{\vee}$ for $1\leq i\leq m$ is constructed
where $V(\tau_i)=\overline{Y_i}$ for $1\leq i\leq m$ are the closures
of the $T$-invariant cells. The piecewise Laurent polynomial
description for this basis seems to agree with that of our basis. More
explicit relations need to be explored.  \erem

\section{Concluding Remarks}

\brem\label{cl1}{\em Concept of cellular varieties:} We observe here
that the concept of cellular varieties is an extension of the concept
of smoothness. Furthermore, we also note that it is a good extension
in the sense that it allows us to extend several nice results on
smooth varieties to the non-smooth setting.  \erem

\brem\label{cl2} {\em Extension to non-complete cellular toric
  varieties:} We can relax the condition of completeness of $\Delta$
to those fans $\Delta$ for which the underlying space of the polytopal
complex $\mathcal{P}_{\Delta}$ is a connected compact topological
manifold with boundary. In particular, $\mathcal{P}_{\Delta}$ is a
pure and strongly connected polytopal complex (see Section \ref{scpc})
such that every ridge lies in either one or two facets. Thus
$\mathcal{P}_{\Delta}$ is an $(n-1)$-dimensional pseudomanifold with
boundary. A fan $\Delta_{\sigma}$ consisting of subdivisions of a
strongly convex rational polyhedral cone in $N$ is one such (see
Example \ref{example1}), since the underlying space of
$\mathcal{P}_{\Delta_{\sigma}}$ is an $(n-1)$-ball. In particular, any
$(n-1)$-dimensional cone in such a fan lies in either one
$n$-dimensional cone or two $n$-dimensional cones. The results of
Section \ref{GKMtoricsection} and Section \ref{pl} can be extended for
such a non-complete fan $\Delta$. The details for this will be
provided elsewhere. It may also be interesting to explore connections
with \cite[Theorem 6.9]{BaBrFiKa}, where a topological criterion is
being discussed for a fan to be equivariantly formal.  \erem

\brem\label{cl3} {\em Extension to non-smooth cellular varieties other
  than toric varieties:} The results in Section \ref{GKMtoricsection}
and Section \ref{basis1} can be generalized to other cellular (not
necessarily toric) varieties. These are work in progress and shall
appear elsewhere. \erem

\brem\label{extension to cellular toric bundles and related spaces}
Recently the results in this paper have been extended to the relative
setting of cellular toric bundles with applications to the description
of the topological equivariant $K$-ring of toroidal horospherical
embeddings (see \cite{u5}). \erem

\noindent {\bf Acknowledgments:} I am very grateful to Prof. Michel
Brion for his valuable comments and suggestions on the earlier
versions of this manuscript. I sincerely thank the anonymous referee
for a very careful reading and several insightful and invaluable
comments and suggestions for improvement of the manuscript. In
particular, the superfluous assumption that the toric variety is
simplical was removed following the suggestions of the referee.

\noindent {\bf Funding:} I thank SERB MATRICS grant no:
MTR/2022/000484 for financial support.

\end{document}